\documentclass[11pt]{amsart} 
\usepackage{verbatim, latexsym, amssymb, amsmath,color}
\usepackage{epsfig} 

\numberwithin{equation}{section}

\def\R{\mathbb R}

\def\N{\mathbb N}
\def\Z{\mathbb Z}

\def\dmn{\mathrm{dmn}}

\def\s{\mathrm{spt}}

\theoremstyle{remark}

\theoremstyle{definition}

\title{Morse index of multiplicity one min-max minimal hypersurfaces}
\author{Fernando C. Marques and Andr\'e Neves}
\address{Institute for Advanced Study and Princeton University \\ Princeton NJ 08544 \\USA}
\email{coda@ias.edu, coda@math.princeton.edu}
\address{University of Chicago \\ Department of Mathematics \\ Chicago IL 60637\\ USA /Imperial College London\\ Huxley Building \\ 180 Queen's Gate \\ London SW7 2RH \\ United Kingdom}
\email{aneves@uchicago.edu, a.neves@imperial.ac.uk}
\thanks{ The first author is partly supported by NSF-DMS-1811840. The second author is partly supported by NSF  DMS-1710846 and EPSRC Programme Grant EP/K00865X/1.}

\begin{document}

\maketitle

\begin{abstract}
In this paper, we prove that the Morse index of a multiplicity one, smooth, min-max minimal hypersurface is generically equal to the dimension of the homology class detected by the families  used in the construction.  This confirms part of the program (\cite{marques-icm},  \cite{marques-neves-cycles}, \cite{marques-neves-index},  \cite{neves-icm}) proposed by the  authors with the goal of  developing a Morse theory for the area functional.
 \end{abstract}

\section{Introduction}

Morse theory concerns the relationship between the structure of the set of critical points of a typical function and the topology of the space on which the function is defined. It was originally devised by Morse to study geodesics (\cite{milnor}), as these are critical points of the length functional defined in the space of paths. We are motivated by the problem of developing a Morse theory for the $n$-dimensional area functional, defined on the space of closed hypersurfaces of some compact manifold $M^{n+1}$. Minimal hypersurfaces are the critical points.

In 1965, Almgren \cite{almgren-varifolds} devised a very  general min-max theory for the area functional that was later improved by Pitts (\cite{pitts}, 1981) in the hypersurface case. Until recently, the main application of Almgren-Pitts theory was the proof that every compact Riemannian manifold $(M^{n+1},g)$, with dimension $3 \leq (n+1) \leq 7$, contains a smooth, embedded, closed minimal hypersurface. If $(n+1)\geq 8$, the minimal hypersurface can be singular in a codimension 7 set by Schoen-Simon \cite{schoen-simon} regularity theory.

In the last few years, we have been able to develop new tools (\cite{marques-neves-index}, \cite{liokumovich-marques-neves}) and find applications of the theory to a variety of problems (\cite{marques-neves-willmore}, \cite{agol-marques-neves}, \cite{marques-neves-infinitely}, \cite{irie-marques-neves}, \cite{marques-neves-song}). We have introduced the idea of using multiparameter sweepouts of mod two flat cycles (\cite{marques-neves-infinitely}), which allowed the use  of topological techniques inspired by Lusternik-Schnirelmann theory. Recently, we have discovered  in \cite{irie-marques-neves} (with Irie) that for generic metrics on $M^{n+1}$, $3 \leq (n+1) \leq 7$, the union of all smooth, embedded, closed minimal hypersurfaces is dense in $M$. This settles the generic case of Yau's Conjecture (\cite{yau1}) from 1982 , about the existence of infinitely many minimal surfaces in any $(M^3,g)$, by proving that a much stronger property, namely density, holds true. This was an application of our Weyl Law for the Volume Spectrum (\cite{liokumovich-marques-neves}, with Liokumovich), which also leads to the existence of an equidistributed sequence of closed minimal hypersurfaces (\cite{marques-neves-song}, with Song). 

The general case of Yau's Conjecture has finally been solved in exciting new work of Song  \cite{song-infinitely-many}.
 In a beautiful paper, Song was able to localize the methods of \cite{marques-neves-infinitely} and proved the existence of infinitely many minimal hypersurfaces  trapped inside a domain bounded by stable hypersurfaces.

In a series of papers (\cite{marques-icm},  \cite{marques-neves-cycles}, \cite{marques-neves-index},  \cite{neves-icm}), the authors proposed a program to obtain a Morse-theoretic description of the set of minimal hypersurfaces in the generic case.
 The  authors conjectured:
 \medskip

\subsection{Morse Index Conjecture}\label{morse.conjecture} {\it For a generic metric $g$ on $M^{n+1}$,  $3\leq (n+1)\leq 7$, there exists a sequence $\{\Sigma_k\}$ of smooth, embedded, two-sided,  closed minimal hypersurfaces such that:
\begin{itemize}
\item[(1)] ${\rm index}(\Sigma_k)=k$,
\item[(2)] $C^{-1}k^\frac{1}{n+1}\leq {\rm area}(\Sigma_k)\leq Ck^\frac{1}{n+1}$ for some $C>0$.
\end{itemize}  
}
\medskip

The  authors proposed a program to prove this conjecture based on three main components: the use of min-max constructions over multiparameter sweepouts to obtain existence results, the characterization of the Morse index of min-max minimal hypersurfaces under the multiplicity one assumption, and a proof of the Multiplicity One Conjecture:

\subsection{Multiplicity One Conjecture}\label{m.o.c} {\it For generic metrics on $M^{n+1}$,  $3\leq (n+1)\leq 7$,  any component of a  closed, minimal hypersurface obtained by min-max methods is  two-sided and has multiplicity one.}
\medskip

The Morse Index Conjecture follows if one can implement the three parts of the program either in  the original Almgren-Pitts setting or in the Allen-Cahn setting, which can be seen as an $\varepsilon$-regularization. The first part of the program was done in the Almgren-Pitts setting by the authors (\cite{marques-neves-infinitely}, \cite{marques-neves-index}), and in the Allen-Cahn setting by  Guaraco \cite{guaraco} and Gaspar-Guaraco \cite{gaspar-guaraco}.

The  Multiplicity One Conjecture \ref{m.o.c}, adapted to the Allen-Cahn setting, was recently proven when $(n+1)=3$ in exciting work by Chodosh and Mantoulidis \cite{chodosh-mantoulidis}. In \cite{chodosh-mantoulidis}, they also finish the Morse index characterization for  multiplicity one Allen-Cahn minimal hypersurfaces (assuming smoothness) in any dimension. In dimension three, they prove remarkable new curvature estimates and strong sheet separation estimates for stable solutions, building on important work of Wang and Wei \cite{wang-wei}. Combined with the existence theory of \cite{gaspar-guaraco}, this gives a new proof of Yau's Conjecture for generic metrics in dimension three. The authors' program has been carried out completely in dimension three (for the Allen-Cahn setting). In particular, the Morse Index Conjecture is true if $(n+1)=3$  (the sublinear growth of the area was proven in Theorem 3.2 of Gaspar-Guaraco \cite{gaspar-guaraco}). We remark that density and equidistribution of minimal hypersurfaces for generic metrics can be also proven via the Allen-Cahn approach (Gaspar-Guaraco, \cite{gaspar-guaraco-weyl}).

We point out that a  multiplicity one property has been proven by Pigati and Rivi\`{e}re \cite{pigati-riviere} in the parametric setting, in which one considers families of maps defined on a two-dimensional surface and taking values in a compact Riemannian manifold of any dimension. This was done by using the viscosity method introduced by Rivi\`{e}re in \cite{riviere} as an alternative to the perturbed functionals of Sacks and Uhlenbeck \cite{sacks-uhlenbeck}. The parametric approach produces immersed (rather than embedded) minimal surfaces with possible branch points. In \cite{pigati-riviere}, Pigati and Rivi\`{e}re prove upper bounds for the Morse index of the minimal surface and discuss the possibility of establishing  lower bounds similar  to those established in our paper.

In this paper, we complete the characterization of the Morse index of Almgren-Pitts min-max minimal hypersurfaces under the multiplicity one assumption. Inspired by finite dimensional Morse theory, one expects that generically the Morse index of a min-max minimal hypersurface $\Sigma$ should be equal to the dimension $k$ of the homology class detected by the families used in the min-max process. The upper bound ${\rm index}(\Sigma)\leq k$ was proven in our previous work \cite{marques-neves-index} for arbitrary metrics, and it was one of the ingredients in the proof of the density result \cite{irie-marques-neves}. In this paper, we prove the lower bound $k\leq {\rm index}(\Sigma)$ for generic (bumpy) metrics under the assumption that the multiplicity of $\Sigma$ is one. This confirms the heuristics coming from finite-dimensional Morse theory, as long as one has the multiplicity one property. Our paper can be seen as a continuation of \cite{marques-neves-index}.

Let us be more precise. Let $(M^{n+1},g)$ be an $(n+1)$-dimensional closed Riemannian manifold,  with $3\leq (n+1)\leq 7$. 
Let $X$ be a simplicial complex of dimension $k$ and $\Phi:X \rightarrow \mathcal Z_n(M^{n+1};{\bf F};\mathbb{Z}_2)$ be a continuous map. Here $\mathcal Z_n(M;\mathbb{Z}_2)$ denotes the space of $n$-dimensional mod 2 flat chains $T$
 in $M$ with $T=\partial U$ for some $(n+1)$-dimensional mod 2 flat chain $U$. These are called {\it flat cycles}. The notation $\mathcal Z_n(M;{\bf F};\mathbb{Z}_2)$ indicates that the space $\mathcal Z_n(M;\mathbb{Z}_2)$ has been endowed with the ${\bf F}$-metric, to be defined later (see Section \ref{preliminaries}). Basically this means $\Phi$ is continuous in both the flat and the varifold topologies.  We denote by $\mathcal{V}_n(M)$ the closure, in the weak topology, of the space of $n$-dimensional rectifiable varifolds with support contained in $M$. A flat cycle $T \in \mathcal Z_n(M;\mathbb{Z}_2)$ induces a varifold $|T|\in \mathcal{V}_n(M)$.
 
We let $\Pi$ be the class of all continuous maps $\Phi':X \rightarrow \mathcal Z_n(M^{n+1};{\bf F};\mathbb{Z}_2)$ such that $\Phi$ and $\Phi'$ are homotopic to each other in the flat topology. 
 The {\it width} of $\Pi$ is defined to be the min-max invariant:
 $$
 {\bf L}(\Pi) = \inf_{\Phi' \in \Pi}\sup_{x\in X}\{{\bf M}(\Phi'(x))\},
 $$
 where ${\bf M}(T)$ denotes the mass of $T$ (or $n$-dimensional area).
 
 Given a sequence  $\{\Phi_i\}$  of continuous maps from $X$ into $\mathcal Z_n(M;{\bf F};\mathbb{Z}_2)$, we set 
$${\bf L}(\{\Phi_i\}):=\limsup_{i \rightarrow \infty} \sup_{x\in X} {\bf M}(\Phi_i(x)).$$
When ${\bf L}(\{\Phi_i\})={\bf L}(\Pi)$, we say $\{\Phi_i\}$ is a {\it minimizing sequence} in $\Pi$. The {\it critical set} of $\{\Phi_i\}$ is the set ${\bf C}(\{\Phi_i\})$  of varifolds $V\in \mathcal{V}_n(M)$ with ${\bf M}(V)={\bf L}(\Pi)$ and such that there exist sequences  $\{i_j\}\subset \{i\}$ and $\{x_j\}\subset X$ with
$$
\lim_{j\rightarrow \infty} {\bf F}(|\Phi_{i_j}(x_j)|,V)=0.
$$

The Almgren-Pitts min-max theory gives that if ${\bf L}(\Pi)>0$ and $3 \leq (n+1) \leq 7$, then there exists a  disjoint collection $\{\Sigma_1, \dots, \Sigma_N\}$ of closed, smooth, embedded, minimal hypersurfaces in $M$ and a set of integers $\{m_1, \dots, m_N\}\subset \mathbb{N}$, such that
$$
V= m_1 \cdot |\Sigma_1| + \cdots + m_N \cdot |\Sigma_N| \in {\bf C}(\{\Phi_i\}).
$$
Any $V\in {\bf C}(\{\Phi_i\}) $ of this form is called an {\it embedded minimal cycle}. The {\it Morse index} of $V$ is the number 
${\rm index}(V)=\sum_{i=1}^N {\rm index}(\Sigma_i)$.  If $m_1=\cdots=m_N=1$, we say $V$ has {\it   multiplicity one}.

It follows indirectly from the work of Almgren \cite{almgren} that $\mathcal Z_n(M;\mathbb{Z}_2)$ is weakly homotopically equivalent to $\mathbb{RP}^\infty$. In Section \ref{space.of.cycles}, we give a simpler and more direct proof of this fact. In particular, $H^1(\mathcal Z_n(M;\Z_2);\Z_2)=\Z_2=\{0,\bar\lambda\}$ and, for every $k\in\N$, we can consider  the set $\mathcal{P}_k$ of all continuous maps $\Phi:Y \rightarrow \mathcal Z_n(M;{\bf F};\Z_2)$  such that $\Phi^*(\bar \lambda)^k$ does not vanish in $H^k(Y;\Z_2)$, where $Y$ is any finite dimensional  simplicial complex. A map in $\mathcal{P}_k$ is called a $k$-sweepout. The {\it $k$-width} of $(M,g)$ is the number 
$$\omega_k(M,g)=\inf_{\Phi\in\mathcal{P}_k}\sup_{x\in dmn(\Phi)}{\bf M}(\Phi(x)),$$
where ${\rm dmn}(\Phi)$ is the domain of $\Phi$ (see Gromov \cite{gromov88}, Guth \cite{guth}, and \cite{marques-neves-infinitely} by the authors).

A Riemannian metric $g$ is  called {\it bumpy} if there is no closed, smooth, immersed, minimal hypersurface that admits a non-trivial Jacobi  field. White showed in \cite{white2, white3} that bumpy metrics are generic in the usual $C^\infty$ Baire sense. In Section \ref{proof.of.main.theorem}, we will prove
that if a metric $g$ is bumpy then for every $k\in \mathbb{N}$ there exists a homotopy class of $k$-sweepouts $\Pi$ such that ${\bf L}(\Pi)=\omega_k(M,g)$.

In this paper, we will prove:

\subsection{Main Theorem}\label{main.theorem}
{\it Suppose $(M^{n+1},g)$ is a bumpy metric, $3\leq (n+1) \leq 7$, and let $\Pi$ be a homotopy class of $k$-sweepouts with ${\bf L}(\Pi)=\omega_k(M,g)$. Suppose $\{\Phi_i\}$ is a minimizing sequence in $\Pi$ such that every embedded minimal cycle of ${\bf C}(\{\Phi_i\})$ has multiplicity one. Then
there exists an embedded minimal cycle  $\Sigma \in {\bf C}(\{\Phi_i\})$ (hence ${\rm area}(\Sigma)=\omega_k(M,g)$)  with
$$
{\rm index}(\Sigma) = k.
$$
}

The main ingredient in the proof is White's Local Min-max Theorem (\cite{white-minmax}). White proved that a nondegenerate, multiplicity one minimal hypersurface of index $k$ is the solution of a local min-max problem over $k$-parameter families with the boundary fixed.  In Section \ref{white.section}, we will prove a version of White's results in which competitors are not necessarily contained in a tubular neighborhood of the original hypersurface (see also Inauen-Marchese \cite{inauen-marchese}).

We further conjecture that the minimal hypersurfaces whose existence is predicted by the Morse Index Conjecture \ref{morse.conjecture} should satisfy:
\begin{equation}\label{weyl.law}
\lim_{k\rightarrow \infty} \frac{{\rm area}(\Sigma_k)}{k^\frac{1}{n+1}} = a(n) {\rm vol}(M,g)^\frac{n}{n+1},
\end{equation}
for some universal constant $a(n)>0$.
Our program, if carried out in the Almgren-Pitts setting, would produce a minimal hypersurface $\Sigma_k$ with ${\rm index}(\Sigma_k)=k$ and ${\rm area}(\Sigma_k)=\omega_k(M,g)$. Hence  identity (\ref{weyl.law}) would follow immediately from the  Weyl Law for the Volume Spectrum, proven by the authors with Liokumovich in \cite{liokumovich-marques-neves}. Property (\ref{weyl.law}) follows in dimension three from \cite{chodosh-mantoulidis}, \cite{gaspar-guaraco} and the Weyl Law proven in Gaspar-Guaraco \cite{gaspar-guaraco-weyl}.

The Multiplicity One Conjecture  \ref{m.o.c} was proven for the one-parameter min-max setting ($k=1$) in \cite{marques-neves-index}, when $M$ has no one-sided embedded hypersurfaces. Previously, the one-parameter case in manifolds with positive Ricci curvature had been studied before by Zhou  \cite{zhou, zhou2}, who did important work extending to high dimensions previous work of  the authors \cite{marques-neves-duke}. These results have been improved recently by Ketover and the authors in \cite{ketover-marques-neves}, where new index characterizations and multiplicity one theorems are proven as a consequence of the catenoid estimate. Related results have been proven for the least-area closed minimal hypersurface by Mazet and Rosenberg \cite{mazet-rosenberg} and later by Song \cite{song}.

General upper bounds for the Morse index of minimal hypersurfaces that are limit interfaces of solutions to the Allen-Cahn equation have been given by Hiesmayr \cite{hiesmayr} and Gaspar \cite{gaspar}. Lower index bounds for the one-parameter min-max case in the Simon-Smith \cite{smith} setting were proven by Ketover-Liokumovich \cite{ketover-liokumovich} and Song \cite{song}. It would be interesting to know if similar index bounds hold in  other min-max constructions (e.g. Chambers-Liokumovich \cite{chambers-liokumovich},  De Lellis-Ramic \cite{delellis-ramic}, De Lellis-Tasnady \cite{delellis-tasnady}, Li-Zhou \cite{li-zhou-free}, Montezuma \cite{montezuma}-\cite{montezuma.mountainpass},  Zhou-Zhu \cite{zhou-zhu}). 

If $(n+1)=2$, min-max theory produces stationary geodesic networks and not closed geodesics (Aiex \cite{aiex}, Mantoulidis \cite{mantoulidis}). In higher codimension, it produces stationary integral varifolds (Almgren \cite{almgren-varifolds}, Stern \cite{stern1}-\cite{stern2}) and the optimal regularity is still an open problem. We conjecture that the min-max minimal variety should be smooth outside a set of codimension two.

 In this paper, we will use Sharp's compactness theorem for minimal hypersurfaces with both area and index uniformly bounded from above (see also Carlotto \cite{carlotto}, Chodosh-Ketover-Maximo \cite{chodosh-ketover-maximo} and Li-Zhou \cite{li-zhou} for other compactness theorems).

The methods of this paper apply equally well, with minor modifications, to situations in which one performs min-max over families (with boundary) that detect some nontrivial relative homology class of the space of cycles, with mod two or integer coefficients.

The paper is organized as follows. In Section \ref{preliminaries}, we introduce the basic notation of Geometric Measure Theory and Almgren-Pitts min-max theory. In Section \ref{interpolation}, we develop further the interpolation techniques introduced in \cite{marques-neves-willmore}. In Section \ref{combinatorial}, we improve the combinatorial argument of Pitts and prove  the existence of minimizing sequences such that every element of the critical set is almost smooth. In Section \ref{space.of.cycles}, we give a direct proof that the space of mod 2 cycles is weakly homotopically equivalent to $\mathbb{RP}^\infty$. In Section \ref{white.section}, we prove an extension of White's Local Min-max Theorem (\cite{white-minmax}). In Section \ref{proof.of.main.theorem}, we prove the Main Theorem \ref{main.theorem}. 

We have added an Addendum to reflect Zhou's recent solution (\cite{zhou-multiplicity}) of the Multiplicity One Conjecture.

\section{Preliminaries}\label{preliminaries}

Let $(M^{n+1},g)$ be an $(n+1)$-dimensional closed Riemannian manifold.  We assume, for convenience, that $(M,g)$ is isometrically embedded in some Euclidean space $\mathbb{R}^J$.

 The spaces we will work with in this paper are:
\begin{itemize}
\item the space ${\bf I}_{l}(M;\mathbb{Z}_2)$  of $l$-dimensional  flat chains    in $\mathbb{R}^J$ with coefficients in $\mathbb{Z}_2$ and support contained  in $M$, where $l=n$ or  $n+1$; 
\item the space ${\mathcal Z}_n(M;\mathbb{Z}_2)$  of flat chains  $T \in {\bf I}_n(M;\mathbb{Z}_2)$ such that there exists $U \in {\bf I}_{n+1}(M;\mathbb{Z}_2)$ with  $\partial U=T$;
\item the closure $\mathcal{V}_n(M)$, in the weak topology, of the space of $n$-dimensional rectifiable varifolds in $\mathbb{R}^J$ with support contained in $M$. 
\end{itemize}
We assume implicitly that ${\bf M}(T)+{\bf M}(\partial T)<\infty$ for every $T\in {\bf I}_{l}(M;\mathbb{Z}_2)$. We will refer to ${\mathcal Z}_n(M;\mathbb{Z}_2)$ as the {\it space of cycles}. Flat chains over a finite coefficient group were introduced by Fleming \cite{fleming}.

Given $T\in {\bf I}_l(M;\mathbb{Z}_2)$,  we denote by $|T|$ and $||T||$ the integral varifold   and the Radon measure in $M$ associated with $|T|$, respectively;  given $V\in \mathcal{V}_n(M)$, $||V||$ denotes the Radon measure in $M$ associated with $V$.  The space of $n$-dimensional integral 
varifolds with support in $M$ is denoted by $\mathcal{IV}_n(M)$.

 The  spaces above come with several relevant metrics. The {\it mass} of $T \in {\bf I}_l(M;\mathbb{Z}_2)$ is denoted by ${\bf M}(T)$, and the metric ${\bf M}(T_1,T_2)={\bf M}(T_1-T_2)$ defines the mass topology. The  {\it flat metric} 
 $$
 \mathcal F(T_1,T_2) = \inf \{{\bf M}(Q)+{\bf M}(R): T_1-T_2=Q+\partial R\}
 $$
 induces the flat topology (we put 
$\mathcal{F}(T)=\mathcal{F}(T,0)$). The  ${\bf F}$-{\it metric}  is defined in  the book of Pitts  \cite[page 66]{pitts} and   induces the varifold weak topology on $\mathcal{V}_n(M)\cap \{V: ||V||(M) \leq a\}$ for any $a$. It satisfies 
$$||V||(M) \leq ||W||(M) +{\bf F}(V,W)$$ for all $V,W \in \mathcal{V}_n(M)$. We denote by ${\overline{\bf B}^{\bf F}_{\delta}(V)}$ and ${\bf B}^{\bf F}_{\delta}(V)$ the closed and open metric balls, respectively, with radius $\delta$ and center $V \in \mathcal{V}_n(M)$. Similarly, we denote by ${\overline{\bf B}^{\mathcal F}_{\delta}(T)}$ and ${\bf B}^{\mathcal F}_{\delta}(T)$ the corresponding balls with  center $T \in \mathcal{Z}_n(M;\mathbb{Z}_2)$ in the flat metric. 
Finally,  the ${\bf F}$-{\it metric} on ${\bf I}_l(M;\mathbb{Z}_2)$ is defined by
$$ {\bf F}(S,T)=\mathcal{F}(S-T)+{\bf F}(|S|,|T|).$$
We have ${\bf F}(|S|,|T|) \leq {\bf M}(S,T)$ and hence $ {\bf F}(S,T) \leq 2{\bf M}(S,T)$ for any $S,T \in {\bf I}_l(M;\mathbb{Z}_2)$.

We assume that  ${\bf I}_l(M;\mathbb{Z}_2)$ and  ${\mathcal Z}_n(M;\mathbb{Z}_2)$ have the topology induced by the flat metric. When endowed with
the topology of the ${\bf F}$-metric or the mass norm, these spaces will be denoted by  ${\bf I}_l(M;{\bf F};\mathbb{Z}_2)$, ${\mathcal Z}_n(M;{\bf F};\mathbb{Z}_2)$, ${\bf I}_l(M;{\bf M};\mathbb{Z}_2)$, ${\mathcal Z}_n(M;{\bf M};\mathbb{Z}_2)$, respectively. The space $\mathcal{V}_n(M)$ is considered with the weak topology of varifolds.
  
  We denote by $B(p,r)$ the Euclidean ball of radius $r$ and center $p$ in $\mathbb{R}^J$, and by $A(p,s,t)$ the annulus $B(p,t)\setminus \bar B(p,s)$.

\subsection{Cubical complexes} 

We denote by  $I^m=[0,1]^m$  the $m$-dimensional unit cube. 
For each $j\in \N$, $I(1,j)$ denotes the cell complex on $I^1$  whose $1$-cells and $0$-cells (also called vertices) are, respectively,  
$$[0,3^{-j}], [3^{-j},2 \cdot 3^{-j}],\ldots,[1-3^{-j}, 1]\quad\mbox{and}\quad [0], [3^{-j}],\ldots,[1-3^{-j}], [1].$$

We consider the $m$-dimensional cubical complex on $I^m$: 
$$I(m,j)=I(1,j)\otimes\ldots \otimes I(1,j)\quad (\mbox{$m$ times}).$$
Then $\alpha=\alpha_1 \otimes \cdots\otimes \alpha_m$ is a $p$-cell of $I(m,j)$ if and only if $\alpha_i$ is a cell
of $I(1,j)$ for each $i$, and $\sum_{i=1}^m {\rm dim}(\alpha_i) =p$. We often abuse notation by identifying  a $p$-cell $\alpha$ with its support: $\alpha_1 \times \cdots \times \alpha_m \subset I^m$. $I(m,j)_p$ denotes the set of all $p$-cells in $I(m,j)$, and $I_0(m,j)$ denotes the set of all cells  of $I(m,j)$ that are contained in the boundary  $\partial I^m$. If $Y$ is a subcomplex of $I(m,j)$, we denote by $Y_0$ the set of vertices of $Y$. We denote by $Y(q)$ the subcomplex of $I(m,j+q)$ with support equal to $Y$. We will also need the map ${\bf n}(i,j):I(m,i)_0\rightarrow I(m,j)_0$, $j\leq i$: ${\bf n}(i,j)(x)$ is the closest vertex in $I(m,j)_0$ to $x$.

Our parameter space $X$ will be a  cubical complex of dimension $k$, meaning  a subcomplex of dimension  $k$ of $I(m,j)$ for some $m$ and $j$. Every such cubical complex  is homeomorphic to a finite simplicial complex and vice-versa (see Chapter 4 of \cite{buchstaber-panov}). 

\subsection{Definition} Let $\Phi:X \rightarrow \mathcal Z_n(M^{n+1};{\bf F};\mathbb{Z}_2)$ be a continuous map. The {\it homotopy class} of $\Phi$ is the class  $\Pi$  of all continuous maps $\Phi':X \rightarrow \mathcal Z_n(M^{n+1};{\bf F};\mathbb{Z}_2)$ such that $\Phi$ and $\Phi'$ are homotopic to each other in the flat topology. 

\subsection{Remark} Notice that our definition of homotopy class is slightly unusual, as we allow homotopies that are continuous in a weaker topology. 

\subsection{Definition} The {\it width} of
$\Pi$ is defined by:
 $$
 {\bf L}(\Pi) = \inf_{\Phi \in \Pi}\sup_{x\in X}\{{\bf M}(\Phi(x))\}.
 $$

\subsection{Definition} A sequence $\{\Phi_i\}\subset \Pi$ is called a {\it minimizing sequence} if 
$$
{\bf L}(\Phi_i):=\sup_{x\in X}{\bf M}(\Phi_i(x))
$$
satisfies ${\bf L}(\{\Phi_i\}):= \limsup_{i\rightarrow \infty} {\bf L}(\Phi_i)={\bf L}(\Pi)$. Any sequence $\{\Phi_{i_j}(x_j)\}$ with $\lim_{j \rightarrow \infty} {\bf M}(\Phi_{i_j}(x_j)) ={\bf L}(\Pi)$, where  $\{i_j\}\subset \{i\}$ is a subsequence and $\{x_j\} \subset X$, is called a {\it min-max sequence}.

\subsection{Definition}\label{image.set} The {\it image set} of $\{\Phi_i\}$ is defined by
\begin{eqnarray*}
{\bf \Lambda}(\{\Phi_i\})&=& \{ V \in \mathcal V_n(M): \exists {\rm \, sequences\, } \{i_j\}\rightarrow \infty, x_{j}\in X\\
&&  \hspace{1cm} {\rm such \, that} \lim_{j\rightarrow \infty} {\bf F}(|\Phi_{i_j}(x_{j})|,V)=0\}.
\end{eqnarray*}

\subsection{Definition}\label{critical.set} If $\{\Phi_i\}$ is a  minimizing sequence in $\Pi$, with $L={\bf L}(\{\Phi_i\})$, the {\it critical set} 
of $\{\Phi_i\}$ is defined by
$$
{\bf C}(\{\Phi_i\})=\{V\in {\bf \Lambda}(\{\Phi_i\}): ||V||(M)=L\}.
$$

\subsection{Pull-tight}\label{pulltight} Following Pitts (\cite{pitts} p.153), we can define  for each $\varepsilon>0$ a continuous map 
\begin{eqnarray*}
H: I \times (\mathcal Z_n(M^{n+1};{\bf F};\mathbb{Z}_2) &\cap& \{T: {\bf M}(T)\leq 2{\bf L}(\Pi)\}) \\
&\rightarrow& \mathcal Z_n(M^{n+1};{\bf F};\mathbb{Z}_2) \cap \{T: {\bf M}(T)\leq 2{\bf L}(\Pi)\}
\end{eqnarray*}
such that:
\begin{itemize}
\item $H(0,T)=T$ for all $T$;
\item $H(t,T)=T$ for all $t\in [0,1]$ if $|T|$ is stationary;
\item ${\bf M}(H(1,T))< {\bf M}(T)$ if $|T|$ is not stationary;
\item ${\bf F}(T, H(t,T))\leq \varepsilon$ for all $t$ and $T$.
\end{itemize}

Given a minimizing sequence $\{\Phi_i^*\}\subset \Pi$, we define $\Phi_i(x)=H(1,\Phi_i^*(x))$ for every $x\in X$. Then $\{\Phi_i\}\subset \Pi$ is also
a minimizing sequence. It follows from the construction that ${\bf C}(\{\Phi_i\}) \subset {\bf C}(\{\Phi_i^*\})$ and that every element of 
${\bf C}(\{\Phi_i\})$ is stationary. 

\subsection{Definition} Any minimizing sequence $\{\Phi_i\}\subset \Pi$ such that every element of ${\bf C}(\{\Phi_i\})$ is stationary is called {\it pulled-tight.}

The next result follows from the Almgren-Pitts (\cite{almgren-varifolds}, \cite{pitts}) min-max theory together with the regularity theory of Schoen-Simon (\cite{schoen-simon}). See Section 3 of \cite{marques-neves-index} for the  formulation presented below.

\subsection{Min-max Theorem}\label{minmax.continuous.thm} {\em Suppose ${\bf L}(\Pi)>0$, and let $\{\Phi_i\}$ be a minimizing sequence in $\Pi$. Then there exists a stationary integral varifold  $V\in {\bf C}(\{\Phi_i\})$ (hence  $||V||(M)={\bf L}(\Pi)$), with support a  closed  minimal hypersurface that is smooth and embedded outside a set of Hausdorff dimension less than or equal to $n-7$.}
\medskip

If $3 \leq (n+1) \leq 7$, we conclude that  there is a disjoint collection $\{\Sigma_1, \dots, \Sigma_N\}$ of closed, smooth, embedded, minimal hypersurfaces in $M$ and a set of integers $\{m_1, \dots, m_N\}\subset \mathbb{N}$, such that
$$
V= m_1 \cdot |\Sigma_1| + \cdots + m_N \cdot |\Sigma_N|.
$$
Any $V\in {\bf C}(\{\Phi_i\}) $ of this form is called an {\it embedded minimal cycle}. The {\it Morse index} of $V$ is the number 
${\rm index}(V)=\sum_{i=1}^N {\rm index}(\Sigma_i)$.  If $m_1=\cdots=m_N=1$, we say $V$ has {\it   multiplicity one}.

\section{Interpolation}\label{interpolation}

We develop further the interpolation machinery initiated in \cite{marques-neves-willmore}. This section is mostly technical and can be skipped in a first reading.

\subsection{Proposition}\label{flat.approximation} {\em Let $X$ be a cubical subcomplex of  $I(m,l)$, and let $\Phi:X\rightarrow \mathcal{Z}_n(M;\mathbb{Z}_2)$ be
a continuous map in the flat topology with no concentration of mass:
$$
\lim_{r\rightarrow 0}\sup_{x\in X}\{{\bf M}(\Phi(x)\llcorner B(p,r)): p\in M\}=0.
$$
There exists a sequence of maps continuous in the mass topology
$$\Phi_i:X \rightarrow \mathcal{Z}_n(M;{\bf M};\mathbb{Z}_2),\quad  i\in\N, $$
that are homotopic to $\Phi$ in the flat topology, and such that 
 $$\sup\{ \mathcal F(\Phi_i(x),\Phi(x)): x\in X\}\rightarrow 0$$
 as $i\rightarrow \infty$.
 }

\begin{proof}
Suppose $\Phi: X \rightarrow \mathcal Z_n(M;\mathbb{Z}_2)$ is as in the statement of the Theorem. Theorem 3.9 of \cite{marques-neves-infinitely} implies there exist a sequence of maps
$$\phi_i:X(k_i)_0 \rightarrow \mathcal{Z}_n(M;\mathbb{Z}_2),$$
with $k_i<k_{i+1}$, and a
sequence of positive numbers $\{\delta_i\}$ converging to zero such that
\begin{itemize}
\item[(i)] $$S=\{\phi_i\}$$ is an $(X,{\bf M})$-homotopy sequence of mappings into $\mathcal{Z}_n(M;{\bf M};\mathbb{Z}_2)$ with ${\bf f}(\phi_i)<\delta_i$ (see Section 2 of \cite{marques-neves-infinitely} for a definition);
\item[(ii)] $$\sup\{\mathcal F(\phi_i(x)-\Phi(x)): x\in X(k_i)_0\}\leq \delta_i;$$
\item[(iii)]$$\sup\{{\bf M}(\phi_i(x)): x\in X(k_i)_0\}\leq \sup\{{\bf M}(\Phi(x)): x\in X\}+\delta_i.$$
\end{itemize}

For sufficiently large $i$, by Theorem 3.10 of \cite{marques-neves-infinitely} we can consider the Almgren extension of $\phi_i$, $\Phi_i: X \rightarrow \mathcal Z_n(M;{\bf M};\mathbb{Z}_2)$. Since $\Phi$ is continuous in the flat metric,  Theorem 3.10 (iii) of \cite{marques-neves-infinitely} implies that
  $$\sup\{ \mathcal F(\Phi_i(x),\Phi(x)): x\in X\}\rightarrow 0$$
  as $i\rightarrow \infty$.
   By Proposition 3.5 of 
 \cite{marques-neves-infinitely}, $\Phi_i$ and $\Phi$ are homotopic to each other in the flat topology.
  
\end{proof}

\subsection{Proposition}\label{close.implies.homotopic} \textit{Let $X$ be a cubical subcomplex of  $I(m,l)$. There exist $\delta_1=\delta_1(M,m)>0$ and $C_1=C_1(M,m)>0$ with the following property:}

\textit{ 
If $\Phi_0,\Phi_1:X \rightarrow \mathcal{Z}_n(M;{\bf M}; \Z_2)$ are continuous maps in the mass topology such that 
$$\sup\{{\bf M}(\Phi_0(x)-\Phi_1(x)):x\in X\} < \delta_1,$$ then there exists a homotopy $H:[0,1] \times X \rightarrow \mathcal{Z}_n(M;{\bf M}; \Z_2)$ with  $H(0,\cdot) = \Phi_0$, $H(1,\cdot)=\Phi_1$ and such that 
\begin{eqnarray*}
&&\sup \{{\bf M}(H(t,x)-\Phi_1(x)): t \in [0,1], x\in X\} \\
&&\hspace{2cm}\leq C_1 \sup\{{\bf M}(\Phi_0(x)-\Phi_1(x)):x\in X\}.
\end{eqnarray*}
 }
\medskip 

\begin{proof}
It follows from combining the work of Almgren (\cite{almgren}, Theorem 8.2) with the mass-continuous deformation map of Pitts (Theorem 4.5 of \cite{pitts}, Section 14 of \cite{marques-neves-willmore}) that there exist  $\delta',C'>0$, depending only on $M$ and $m$, such that if $\Phi:I^k \rightarrow \mathcal{Z}_n(M;\Z_2)$, $k \leq m$, is continuous in the mass topology, $\Phi(x)=0$ for all $x \in \partial I^k$ and ${\bf M}(\Phi(x))\leq\delta'$ for every $x\in I^k$, then there exists a homotopy $H:I^{k+1}\rightarrow \mathcal{Z}_n(M;\Z_2)$ with the following properties:
\begin{itemize}
\item $H$ is continuous in the mass topology;
\item $H(x,0)=0$ and $H(x,1)=\Phi(x)$ for every $x \in I^k$;
\item $H(x,t)=0$ for every $x\in \partial I^k$ and $t\in [0,1]$;
\item $\sup\{{\bf M}(H(w)):w\in I^{k+1}\} \leq C'\sup\{{\bf M}(\Phi(x)):x\in I^{k}\}$.
\end{itemize}
The proof continues as in the Appendix A of \cite{marques-neves-infinitely}.
\end{proof}

If $\phi$ is any map taking values in $\mathcal{Z}_n(M;\mathbb{Z}_2)$ with domain ${\rm dmn}(\phi)$, we define
$$
{\bf m}(\phi,r) = \sup_{x\in {\rm dmn}(\phi)} \{{\bf M}(\phi(x) \llcorner B(p,r)):p\in M\}.
$$

Let $a(k)=2^{-4(k+2)^2-2}$.
The following proposition was proved in \cite{marques-neves-willmore} (Proposition 13.5) when $n=2$ and for integer coefficients, but the same proof applies for $\mathbb{Z}_2$ coefficients in any dimension.
\subsection{Proposition}\label{singleT}  \textit{Let  $\delta,r,P$ be positive real numbers, and  let
 $$T\in  \mathcal{Z}_n(M;\mathbb{Z}_2)\cap\{S:{\bf M}(S)\leq 2P\}.$$  There exist $0<\varepsilon=\varepsilon(T,\delta,r,P)<\delta$,  $q=q(T,\delta,r,P)\in \N$, and $b(k)>0$ for which  the following holds:
given $0<s<\varepsilon$,   $j,l\in \N$ with  $l\leq k+1,$ and 
$$\phi:I_0(l,j)_0\rightarrow  {\bf B}_s^{\mathcal{F}}(T)\cap\{S:{\bf M}(S)\leq 2P-\delta\}$$
 with $$2^{k+2}({\bf m}(\phi,r)+a(k)\delta)\leq \delta/4,$$
there exists
$$\tilde\phi:I(l,j+q)_0\rightarrow  {\bf B}_{s}^{\mathcal{F}}(T)$$
with 
\begin{itemize}
\item[(i)]  ${\bf f}(\tilde\phi) \leq \delta$ if $l=1$ and ${\bf f}(\tilde \phi)\leq b(k)({\bf f}(\phi)+\delta)$ if $l\neq 1$; 
\item[(ii)] $\tilde\phi=\phi\circ {\bf n}(j+q,j)$ on $I_0(l,j+q)_0$;
\item[(iii)] $$\sup_{x\in I(l,j+q)_0}\{{\bf M}(\tilde\phi(x))\}\leq \sup_{x\in I_0(l,j)_0}\{{\bf M}(\phi(x))\}+\delta;$$
\item[(iv)]${\bf m}(\tilde\phi,r)\leq 2^{k+2}({\bf m}(\phi,r)+a(k)\delta).$
\end{itemize}
}

\medskip

Proposition \ref{singleT} has a counterpart in the continuous setting.

\subsection{Proposition}\label{extension.mass}\textit{Let $T\in \mathcal{Z}_n(M;\mathbb{Z}_2)$, and $l\in \mathbb{N}$. Given $\gamma>0$, there exists $0<\eta<\gamma$ such that for every continuous map $\Phi:\partial I^l \rightarrow \mathcal{Z}_n(M;{\bf M};\mathbb{Z}_2)$ with
$$
\sup_{y\in \partial I^l} {\bf F}(\Phi(y),T)<\eta,
$$
there is an extension $\hat{\Phi}:I^l \rightarrow \mathcal{Z}_n(M;{\bf M};\mathbb{Z}_2)$ of $\Phi$ with
$$
\sup_{x\in I^l} {\bf F}(\hat\Phi(x),T)<\gamma.
$$
}

\begin{proof}
Let $P={\bf M}(T)$ and $k=l-1$. Let $C_0=C_0(M,l)$ and $\delta_0=\delta_0(M,l)$ be the positive constants of Theorem 3.10 of \cite{marques-neves-infinitely}, concerning the existence of Almgren extensions. Let also $C_1=C_1(M,l)$ and $\delta_1=\delta_1(M,l)$ be as in Proposition \ref{close.implies.homotopic}.

Choose $\delta <\min\{P/2, \delta_0/(2+4b(k))\}$ and sufficiently small (depending on $T$ and $\gamma$) so that $2(\delta + C_0 (1+2b(k))\delta)\leq \delta_1$, $2C_1( \delta + C_0 (1+2b(k))\delta) \leq \gamma/4$ and if $S\in \mathcal{Z}_n(M;\mathbb{Z}_2)$ satisfies
\begin{eqnarray*}
&&\mathcal{F}(S,T)\leq  2C_0 (1+2b(k))\delta, \\
&& {\bf M}(S)\leq {\bf M}(T)+ \delta + 2C_0 (1+2b(k))\delta,
\end{eqnarray*}
then ${\bf F}(S,T)\leq \gamma/2$. The existence of such  a $\delta$ follows from the fact that convergence in the ${\bf F}$ metric is equivalent to convergence in the flat metric plus convergence of the masses (see 2.1.20 of \cite{pitts}).
Choose $r$ so that $$2^{k+2}({\bf m}(T,2r)+ a(k) \delta) \leq \delta/8.$$

Choose $\eta>0$ sufficiently small so that
\begin{itemize}
\item $\eta < \min\{\gamma, \varepsilon=\varepsilon(T,\delta,r,P), P/2, C_0 (1+2b(k))\delta\}$,
\item $2C_1( \delta + C_0 (1+2b(k))\delta)+\eta \leq \gamma/2$,
\item $2^{k+2}({\bf m}(S,r)+ a(k) \delta) \leq \delta/4$ for any $S \in {\bf B}^{\bf F}_\eta(T)$.
\end{itemize}

Let $\Phi:\partial I^l \rightarrow \mathcal{Z}_n(M;{\bf M};\mathbb{Z}_2)$ be a continuous map that satisfies
$$
\sup_{y\in \partial I^l} {\bf F}(\Phi(y),T)<\eta.
$$
We have
$$\sup_{y \in \partial I^l}{\bf M}(\Phi(y))\leq {\bf M}(T)+\eta\leq 2P-\delta.$$

Choose $j$ sufficiently large so that ${\bf M}(\Phi(y_1)-\Phi(y_2))<\delta$ whenever $y_1,y_2\in \partial I^l$ belong to a common cell of $I_0(l,j)$. Define $\phi=\Phi_{|I_0(l,j)_0}$, hence ${\bf f}(\phi)<\delta$.

We apply Proposition \ref{singleT} to $\phi$, with $s=\eta$, to obtain a map $$\tilde{\phi}: I(l,j+q)_0\rightarrow {\bf B}^\mathcal{F}_\eta(T),$$
$q=q(T,\delta,r, P)$, that satisfies
\begin{itemize}
\item[(i)] ${\bf f}(\tilde{\phi})\leq (1+2b(k))\delta\leq \delta_0/2$;
\item[(ii)] $\tilde\phi=\phi\circ {\bf n}(j+q,j)$ on $I_0(l,j+q)_0$;
\item[(iii)] $$\sup_{x\in I(l,j+q)_0}\{{\bf M}(\tilde\phi(x))\}\leq  \sup_{y\in \partial I^l} {\bf M}(\Phi(y))+\delta.$$
\end{itemize}

Theorem 3.10 of \cite{marques-neves-infinitely} applied to $\tilde{\phi}$ gives a continuous map $\tilde{\Phi}: I^l \rightarrow \mathcal{Z}_n(M;{\bf M};\mathbb{Z}_2)$ such that $\tilde\Phi(x)=\tilde\phi(x)$ for every $x\in I(l,j+q)_0$ and
$$ {\bf M}(\tilde\Phi(x_1)-\tilde\Phi(x_2))\leq C_0{\bf f}(\tilde\phi)\leq  C_0 (1+2b(k))\delta$$ if $x_1,x_2$  lie in a common cell in $I(l,j+q)$.

Given $x\in I^l$, we can choose $\tilde{x}\in I(l,j+q)_0$ belonging  to a cell that contains $x$. 
Then
\begin{eqnarray*}
\mathcal{F}(\tilde\Phi(x),T)&\leq& \mathcal{F}(\tilde\Phi(\tilde x),T)+ \mathcal{F}(\tilde\Phi(\tilde x), \tilde\Phi(x))\\
&\leq& 
\mathcal{F}(\tilde\phi(\tilde x),T)+ {\bf M}(\tilde\Phi(\tilde x), \tilde\Phi(x)) \\
&\leq&  \eta + C_0 (1+2b(k))\delta\\
&\leq&  2C_0 (1+2b(k))\delta.
\end{eqnarray*}
Also,
\begin{eqnarray*}
{\bf M}(\tilde\Phi(x))&\leq& {\bf M}(\tilde\Phi(\tilde x))+{\bf M}(\tilde\Phi(\tilde x), \tilde\Phi(x))\\
&\leq& {\bf M}(\tilde\phi(\tilde x))+C_0 (1+2b(k))\delta\\
&\leq& \sup_{y \in \partial I^l}{\bf M}(\Phi(y))+ \delta + C_0 (1+2b(k))\delta\\
&\leq& {\bf M}(T)+\eta+ \delta + C_0 (1+2b(k))\delta\\
&\leq& {\bf M}(T)+\delta + 2C_0 (1+2b(k))\delta.
\end{eqnarray*}
Then
${\bf F}(\tilde\Phi(x),T) \leq \gamma/2$ for every $x\in I^l$.

For $y \in \partial I^l$, we can choose $\tilde{y}\in I_0(l,j+q)_0$ belonging to a cell in $I(l,j+q)$ containing $y$ and $\tilde{\tilde{y}}={\bf n}(j+q,j)(\tilde{y})\in I_0(l,j)_0$ belonging to a cell in $I(l,j)$ containing $y$. Then
\begin{eqnarray*}
{\bf M}(\Phi(y)-\tilde\Phi(y))&\leq& {\bf M}(\Phi(y)-\tilde\Phi(\tilde y))+ C_0 (1+2b(k))\delta\\
&=& {\bf M}(\Phi(y)-\tilde\phi(\tilde y))+ C_0 (1+2b(k))\delta\\
&=&{\bf M}(\Phi(y)-\phi(\tilde{\tilde{y}}))+C_0 (1+2b(k))\delta\\
&=&{\bf M}(\Phi(y)-\Phi(\tilde{\tilde{y}}))+ C_0 (1+2b(k))\delta\\
&\leq& \delta + C_0 (1+2b(k))\delta,
\end{eqnarray*}
hence ${\bf M}(\Phi(y)-\tilde\Phi(y)) \leq \delta_1/2$.
By Proposition \ref{close.implies.homotopic}, there exists a homotopy $H:[0,1] \times \partial I^l \rightarrow \mathcal{Z}_n(M;{\bf M}; \Z_2)$, with  $H(0,\cdot) = \tilde\Phi$ and $H(1,\cdot)=\Phi$, such that 
\begin{eqnarray*}
&&\sup \{{\bf M}(H(t,y)-\Phi(y)): t \in [0,1], y\in \partial I^l\} \\
&&\hspace{2cm}\leq C_1 \sup\{{\bf M}(\Phi(y)-\tilde\Phi(y)):y\in \partial I^l\} \leq C_1( \delta + C_0 (1+2b(k))\delta).
\end{eqnarray*}
Now
\begin{eqnarray*}
{\bf F}(H(t,y), T) &\leq& {\bf F}(H(t,y), \Phi(y)) + {\bf F}(\Phi(y),T)\\
&\leq& 2{\bf M}(H(t,y), \Phi(y))+ \eta\\
&\leq& 2C_1( \delta + C_0 (1+2b(k))\delta)+\eta\\
&\leq& \gamma/2.
\end{eqnarray*}

Now we can define $\hat{\Phi}:I^l \rightarrow \mathcal{Z}_n(M;{\bf M}; \Z_2)$ by
\begin{itemize}
\item $\hat{\Phi}(x)=\tilde{\Phi}(c+2(x-c))$ if $x \in c+ [-1/4,1/4]^l$,
\item  $\hat{\Phi}(x)=H(\Lambda(x))$ if  $x \notin c+ [-1/4,1/4]^l$,
\end{itemize}
 where $c=(1/2, \dots, 1/2) \in I^l$ and $\Lambda: I^l \setminus (c+ [-1/4,1/4]^l) \rightarrow [0,1] \times \partial I^l$ is any homeomorphism with 
\begin{eqnarray*}
\Lambda(x) = (0, c+2(x-c))
\end{eqnarray*}
for every $x\in \partial(c+ [-1/4,1/4]^l)$,
and 
\begin{eqnarray*}
\Lambda(x)=(1,x)
\end{eqnarray*}
for every $x \in \partial I^l$. The map $\hat{\Phi}$ is our extension.

\end{proof}

\subsection{Proposition}\label{homotopy.F.ball} \textit{Let $T\in \mathcal{Z}_n(M;\mathbb{Z}_2)$, and $k\in \mathbb{N}$. Given $\gamma>0$, there exists $\eta>0$ such that for every cubical complex $X$ of dimension at most $k$ and  a pair of continuous maps $\Phi:X \rightarrow \mathcal{Z}_n(M;{\bf M};\mathbb{Z}_2)$, $\Psi:X \rightarrow \mathcal{Z}_n(M;{\bf M};\mathbb{Z}_2)$ with
$$
\sup_{x\in X} {\bf F}(\Phi(x),T)<\eta, \hspace{0.5cm} \sup_{x\in X} {\bf F}(\Psi(x),T)<\eta,
$$
there is a homotopy $H:[0,1] \times X \rightarrow \mathcal{Z}_n(M;{\bf M};\mathbb{Z}_2)$ with $H(0, \cdot)=\Phi$, $H(1, \cdot)=\Psi$ and such that
$$
\sup_{(t,x)\in [0,1] \times X} {\bf F}(H(t,x),T)<\gamma.
$$
}

\begin{proof}
Set $\gamma_{k+1}=\gamma$, and let $\gamma_k=\eta(T,l=k+1,\gamma)$ be the constant given by Proposition \ref{extension.mass}.
Define by backwards induction $\gamma_i=\eta(T,i+1,\gamma_{i+1})$. Then $\eta=\gamma_0$. Let $X^{(i)}$ be the $i$-dimensional skeleton of $X$. We use Proposition \ref{extension.mass} to first define $H$ on $[0,1]\times X^{(0)}$ and then inductively to define
$H$ on $[0,1]\times X^{(i+1)}$ from $H_{|[0,1]\times X^{(i)}}$.
\end{proof}

\subsection{Proposition}\label{homotopy.F.close} \textit{Let $\mathcal{K} \subset \mathcal{Z}_n(M;{\bf F};\mathbb{Z}_2)$ be a compact set, and $k \in \mathbb{N}$. For every $\alpha>0$, there exists $0<\beta<\alpha$ (depending on $\mathcal{K}$) such that for every cubical complex $X$ of dimension at most $k$ and every pair of continuous maps 
$ \Phi:X\rightarrow \mathcal{Z}_n(M;{\bf M};\Z_2)$ and $\Psi:X\rightarrow \mathcal{Z}_n(M;{\bf M};\Z_2)$ with $\Psi(X) \subset {\bf B}^{\bf F}_\beta(\mathcal{K})$ and 
$$\sup\{{\bf F}(\Phi(x),\Psi(x)) : x\in X\}<\beta,$$
there exists a homotopy $H:[0,1]\times X \rightarrow \mathcal{Z}_n(M;{\bf M};\Z_2)$ satisfying  $H(0,\cdot) = \Phi$ and $H(1,\cdot)=\Psi$  and such that
$$
{\bf F}(H(t,x), \Psi(x))\leq \alpha \quad\mbox{for all }(t,x)\in  [0,1]\times X.
$$
}

\begin{proof}
Let us use induction in $k$. If $k=0$, $X$ is a finite set. For every $T\in \mathcal{K}$, let $\eta(T)=\eta(T,l=1,\gamma=\alpha/3)<\alpha/3$ be the constant given by Proposition \ref{extension.mass}. By compactness, there are $\{T_1, \dots, T_q\}\subset \mathcal{K}$ such that
$$\mathcal{K} \subset {\bf B}^{\bf F}_{\eta(T_1)/2}(T_1) \cup \cdots \cup {\bf B}^{\bf F}_{\eta(T_q)/2}(T_q).$$

Let  $\eta=\min\{\eta(T_1), \dots, \eta(T_q)\}$, and choose $\beta=\eta/8$.
Given $x\in X$, there exists $T_i$ such that $\Psi(x)\in {\bf B}^{\bf F}_{\beta+ \eta(T_i)/2}(T_i)$. Then $\Phi(x) \in {\bf B}^{\bf F}_{2\beta+ \eta(T_i)/2}(T_i)$. Since $\{\Phi(x), \Psi(x)\} \subset {\bf B}^{\bf F}_{\eta(T_i)}(T_i)$,  Proposition \ref{extension.mass} gives  a continuous map $H:[0,1] \rightarrow  \mathcal{Z}_n(M;{\bf M};\Z_2)$ with $H(0)=\Phi(x)$, $H(1)=\Psi(x)$ and $H(t) \in {\bf B}^{\bf F}_{\alpha/3}(T_i)$ for every $t \in [0,1]$. Then
$$
{\bf F}(H(t),\Psi(x))\leq {\bf F}(H(t),T_i)+{\bf F}(T_i,\Psi(x))\leq \alpha/3+\eta(T_i)<\alpha
$$
for every $t \in [0,1]$ and we are done.

Now suppose that Proposition \ref{homotopy.F.close} holds for every $k\leq l$. For every $T\in \mathcal{K}$, let $\eta(T)=\eta(T,l+2,\gamma=\alpha/3)<\alpha/3$ be the constant given by Proposition \ref{extension.mass}. By compactness, there are $\{T_1, \dots, T_q\}\subset \mathcal{K}$ such that
$$\mathcal{K} \subset {\bf B}^{\bf F}_{\eta(T_1)/2}(T_1) \cup \cdots \cup {\bf B}^{\bf F}_{\eta(T_q)/2}(T_q).$$ Let  $\eta=\min\{\eta(T_1), \dots, \eta(T_q)\}$. Define $\tilde\alpha=\eta/8$ and let $\tilde{\beta}=\beta(\mathcal{K},l,\tilde\alpha)$ be the constant of Proposition \ref{homotopy.F.close} obtained in the induction step.   Choose $\beta=\min\{\eta/8, \tilde \beta\}.$

Let $X$ be a cubical complex
of dimension $l+1$ and consider  a pair of continuous maps 
$ \Phi:X\rightarrow \mathcal{Z}_n(M;{\bf M};\Z_2)$ and $\Psi:X\rightarrow \mathcal{Z}_n(M;{\bf M};\Z_2)$ with $\Psi(X) \subset {\bf B}^{\bf F}_{\beta}(\mathcal{K})$ and 
$$\sup\{{\bf F}(\Phi(x),\Psi(x)) : x\in X\}<\beta.$$ 
If $j$ is sufficiently large, and $X(j)$ denotes the $j$-th subdivision of $X$, we will have that
$$
{\bf M}(\Phi(x_1),\Phi(x_2))+ {\bf M}(\Psi(x_1),\Psi(x_2)) \leq \eta/4
$$
whenever $x_1$ and $x_2$ belong to a shared cell in $X(j)$.
 
 If $X(j)^{(l)}$ denotes the $l$-dimensional skeleton of $X(j)$, we can apply Proposition \ref{homotopy.F.close} to the restrictions
 $\Phi_{|X(j)^{(l)}}, \Psi_{|X(j)^{(l)}}$. This gives a homotopy $\tilde H:[0,1]\times X(j)^{(l)} \rightarrow \mathcal{Z}_n(M;{\bf M};\Z_2)$ satisfying  $\tilde H(0,\cdot) = \Phi$ and $\tilde H(1,\cdot)=\Psi$  and such that
$$
{\bf F}(\tilde H(t,x), \Psi(x))\leq \tilde \alpha \quad\mbox{for all }(t,x)\in  [0,1]\times X(j)^{(l)}.
$$

Let $\sigma$ be an $(l+1)$-dimensional cell of $X(j)$. We define a continuous map 
$$\Lambda: \partial([0,1]\times \sigma) \rightarrow \mathcal{Z}_n(M;{\bf M};\Z_2)$$
by
\begin{itemize}
\item $\Lambda(0,x) = \Phi(x)$ for every $x\in \sigma$,
\item $\Lambda(1,x)=\Psi(x)$ for every $x\in \sigma$,
\item $\Lambda(t,y)=\tilde H(t,y)$ for every $y\in \partial \sigma$.
\end{itemize}
Choose $x_0\in \sigma$, and let $1\leq i\leq q$ such that $\Psi(x_0)\in {\bf B}^{\bf F}_{\beta+\eta(T_i)/2}(T_i)$. Then $$\Psi(\sigma) \subset {\bf B}^{\bf F}_{\eta/4+\beta+\eta(T_i)/2}(T_i),$$
$$
\Phi(\sigma) \subset {\bf B}^{\bf F}_{\eta/4+2\beta+\eta(T_i)/2}(T_i),
$$
and
 $$\tilde H([0,1]\times \partial \sigma) \subset {\bf B}^{\bf F}_{\tilde\alpha+\eta/4+\beta+\eta(T_i)/2}(T_i).$$
 Hence
 $$
 \Lambda(\partial ([0,1]\times \sigma)) \subset {\bf B}^{\bf F}_{\eta(T_i)}(T_i).
 $$
 Proposition \ref{extension.mass} gives a continuous  extension $\Lambda_\sigma: [0,1]\times \sigma \rightarrow  \mathcal{Z}_n(M;{\bf M};\Z_2)$ of $\Lambda$ such that
 $$
 \Lambda_\sigma([0,1]\times \sigma) \subset {\bf B}^{\bf F}_{\alpha/3}(T_i).
 $$
 Hence ${\bf F}(\Lambda_\sigma(t,x)), \Psi(x)) \leq {\bf F}(\Lambda_\sigma(t,x)), T_i) + {\bf F}(T_i, \Psi(x)) \leq \alpha/3+\eta(T_i)<\alpha.$
 
 We put $H_{|[0,1]\times \sigma}=\Lambda_\sigma$, for every $\sigma$, and we are done.

\end{proof}

The next proposition shows that an ${\bf F}$-continuous map can be arbitrarily approximated by maps that are continuous in the mass topology.

\subsection{Proposition}\label{M.regularization} \textit{Let $X$ be a cubical complex of dimension $k$, and $\Phi:X \rightarrow \mathcal{Z}_n(M;{\bf F};\mathbb{Z}_2)$ be a continuous map. Then for every $\varepsilon>0$ there exists a continuous map $\Phi': X \rightarrow \mathcal{Z}_n(M;{\bf M};\mathbb{Z}_2)$ and a homotopy $H:[0,1] \times X \rightarrow \mathcal{Z}_n(M;{\bf F};\mathbb{Z}_2)$ with $H(0,\cdot)=\Phi'$, $H(1,\cdot)=\Phi$ and such that
$$
\sup_{(t,x)\in [0,1] \times X} {\bf F}(H(t,x),\Phi(x))<\varepsilon.
$$
}

\begin{proof}
By Proposition A.5 of \cite{marques-neves-index}, there exists a sequence of continuous maps $\Phi_i:X \rightarrow \mathcal{Z}_n(M;{\bf M};\mathbb{Z}_2)$ such that
$$
\lim_{i\rightarrow \infty} \sup_{x\in X} {\bf F}(\Phi(x), \Phi_i(x)) =0. 
$$

Let $\alpha>0$ and $\mathcal{K}=\Phi(X)$. Let $0<\beta<\alpha$ be the constant given by Proposition \ref{homotopy.F.close} with $k={\rm dim}(X)$. If $i$ is sufficiently large, then  $\Phi_i(X)\subset {\bf B}_\beta^{\bf F}(\mathcal{K})$
and
$$
\sup_{x\in X} {\bf F}(\Phi_i(x), \Phi_{i+1}(x)) < \beta.
$$
Proposition \ref{homotopy.F.close} gives a homotopy $\hat H:[0,1]\times X \rightarrow \mathcal{Z}_n(M;{\bf M};\Z_2)$ satisfying  $\hat H(0,\cdot) = \Phi_i$ and $\hat H(1,\cdot)=\Phi_{i+1}$  and such that
$$
{\bf F}(\hat H(t,x), \Phi_{i+1}(x))\leq \alpha \quad\mbox{for all }(t,x)\in  [0,1]\times X.
$$

Since $\alpha$ is arbitrary, we can find for every sufficiently large $i$  a homotopy $H_i$ (continuous in the mass norm) between $\Phi_i$ and $\Phi_{i+1}$ such that
$$
\sup_{x\in X} {\bf F}(H_i(t,x),\Phi(x))=\varepsilon_i \rightarrow 0
$$
as $i\rightarrow \infty$. The desired $H$ is obtained by concatenating the sequence $\{H_i\}.$ 

\end{proof}

We now prove an analogous statement to Proposition \ref{homotopy.F.close} for maps that are continuous in the ${\bf F}$-metric.

\subsection{Theorem}\label{homotopy.F.thm} \textit{Let $\mathcal{K} \subset \mathcal{Z}_n(M;{\bf F};\mathbb{Z}_2)$ be a compact set, and $k \in \mathbb{N}$. For every $\alpha>0$, there exists $0<\beta<\alpha$ (depending on $\mathcal{K}$) such that for every cubical complex $X$ of dimension at most $k$ and every pair of continuous maps 
$ \Phi:X\rightarrow \mathcal{Z}_n(M;{\bf F};\Z_2)$ and $\Psi:X\rightarrow \mathcal{Z}_n(M;{\bf F};\Z_2)$ with $\Psi(X) \subset {\bf B}^{\bf F}_\beta(\mathcal{K})$ and 
$$\sup\{{\bf F}(\Phi(x),\Psi(x)) : x\in X\}<\beta,$$
there exists a homotopy $H:[0,1]\times X \rightarrow \mathcal{Z}_n(M;{\bf F};\Z_2)$ satisfying  $H(0,\cdot) = \Phi$ and $H(1,\cdot)=\Psi$  and such that
$$
{\bf F}(H(t,x), \Psi(x))\leq \alpha \quad\mbox{for all }(t,x)\in  [0,1]\times X.
$$
}

\begin{proof}
Given $\alpha>0$,  let $0<\beta=\beta(\mathcal{K}, k, \alpha/2)<\alpha/2$ be the constant given by Proposition \ref{homotopy.F.close}. Let $ \Phi:X\rightarrow \mathcal{Z}_n(M;{\bf F};\Z_2)$ and $\Psi:X\rightarrow \mathcal{Z}_n(M;{\bf F};\Z_2)$ satisfy $\Psi(X) \subset {\bf B}^{\bf F}_\beta(\mathcal{K})$ and 
$$\sup\{{\bf F}(\Phi(x),\Psi(x)) : x\in X\}<\beta.$$

By Proposition \ref{M.regularization}, there exist continuous maps $ \Phi':X\rightarrow \mathcal{Z}_n(M;{\bf M};\Z_2)$ and $\Psi':X\rightarrow \mathcal{Z}_n(M;{\bf M};\Z_2)$ with $\Psi'(X) \subset {\bf B}^{\bf F}_\beta(\mathcal{K})$ and 
\begin{eqnarray*}
&&\sup\{{\bf F}(\Phi'(x),\Psi'(x)) : x\in X\}<\beta,\\
&&\sup\{{\bf F}(\Psi'(x),\Psi(x)) : x\in X\}<\beta,
\end{eqnarray*}
and ${\bf F}$-continuous homotopies $H_1$ between $\Phi'$ and $\Phi$, $H_2$ between $\Psi'$ and $\Psi$ with
\begin{eqnarray*}
{\bf F}(H_1(t,x), \Psi(x))< \beta {\rm \, \, and\, \, } {\bf F}(H_2(t,x), \Psi(x))< \beta \quad\mbox{for all }(t,x)\in  [0,1]\times X.
\end{eqnarray*}

Proposition \ref{homotopy.F.close} then gives an ${\bf M}$-continuous homotopy $H_3$ between $\Phi'$ and $\Psi'$ such that
$$
{\bf F}(H_3(t,x),\Psi'(x)) \leq \alpha/2 \quad\mbox{for all }(t,x)\in  [0,1]\times X.
$$

The desired $H$ is obtained by concatenating $H_1$, $H_2$ and $H_3$.
\end{proof}

\section{Deforming away nonsmooth critical varifolds}\label{combinatorial}

In this section we will prove  the existence of minimizing sequences such that every element of the critical set is almost smooth.

We need to discuss first almost minimizing varifolds, a notion that lies  at the heart of Almgren-Pitts min-max theory.

\subsection{Definition} For each pair of positive numbers $\varepsilon, \delta$, an open set $U\subset \R^J$ and $T\in \mathcal Z_{n}(M, M\setminus U; \Z_2)$, an {\it $(\varepsilon, \delta)$-deformation} of $T$ in $U$ is a finite sequence $\{T_i\}_{i=0}^q$ in $\mathcal Z_{n}(M, M\setminus U;\Z_2)$ with 
\begin{itemize}
\item $T_0=T$ and $\mbox{spt}(T-T_i)\subset U$ for all $i=1,\ldots, q;$
\item ${\bf M}(T_i-T_{i-1})\leq \delta$ for all $i=1,\ldots, q;$
\item ${\bf M}(T_i)\leq {\bf M}(T)+\delta$ for all $i=1,\ldots, q;$
\end{itemize}
and $||T_q||(U)< ||T||(U)-\varepsilon$.

We define $$\mathfrak a (U;\varepsilon,\delta)$$ to be the set of all flat chains  $T\in \mathcal Z_{n}(M, M\setminus U; \Z_2)$ that do not admit $(\varepsilon,\delta)$-deformations in $U$. Loosely speaking this says that every  deformation of  $T$ that is supported in $U$ and that decreases the area by more than $\varepsilon$ must pass through a stage where the area has increased by more than  $\delta$. We have  $\mathfrak a (U;\varepsilon,\delta)\subset \mathfrak a (U;\varepsilon',\delta')$ if $\delta'\leq\delta$ and $\varepsilon\leq \varepsilon'.$

\subsection{Definition} A varifold $V\in \mathcal V_{n}(M)$ is {\em $\Z_2$ almost minimizing} in an open set $U\subset M$ if  for every $\varepsilon>0$, we can find $\delta>0$ and  
$$T\in \mathfrak a (U;\varepsilon,\delta)$$ with ${\bf F}_U(V,|T|)<\varepsilon$.

\medskip

\subsection{Definition} A varifold $V \in \mathcal{V}_n(M)$ is {\it $\mathbb{Z}_2$ almost minimizing in annuli} if 
for each $p \in M$, there exists $r=r(p)>0$ such that $V$ is $\mathbb{Z}_2$ almost minimizing in $M \cap A (p,s,r)$ for all $0<s<r$. 

\medskip

If $V \in \mathcal{V}_n(M)$ is stationary in $M$ and $\mathbb{Z}_2$-almost minimizing in annuli, then $V \in \mathcal{IV}_n(M)$ by Theorem 3.13 of \cite{pitts}. The next result  follows from the  regularity theory of Pitts \cite[Theorem 3.13 and Section 7]{pitts} and Schoen-Simon \cite[Theorem 4]{schoen-simon} (see Theorem 2.11 of \cite{marques-neves-infinitely}).

\subsection{Theorem}\label{regularity.thm} {\em  Let $V  \in \mathcal{V}_n(M)$ be a  stationary varifold that is $\mathbb{Z}_2$ almost minimizing in annuli. Then $V$ is integral and the support of $V$ is a closed,  minimal hypersurface that is smooth and embedded outside a set of codimension 7. } 
\medskip

We now make some comments regarding the existence of varifolds that are almost minimizing in annuli. We are going to use the notation of Section 2 of \cite{marques-neves-infinitely}. This corresponds to the original formulation of Almgren-Pitts min-max theory (\cite{pitts}), which uses discrete families of cycles instead of continuous ones.

Suppose $X$ is a $k$-dimensional cubical subcomplex of $I(m,j)$. 
Consider an $(X,{\bf M})$-homotopy sequence of mappings $S=\{\phi_i\}$.
We set $L={\bf L}(S)$ and assume  that every element in ${\bf C}(S)$ is stationary. 
The basic approach to the existence of almost minimizing varifolds consists in assuming that no element in ${\bf C}(S)$ is almost minimizing in annuli, and  then producing another sequence $S^*=\{\phi^*_i\}$ homotopic to $S$  with the property that ${\bf L}(S^*)<L$. This uses  a combinatorial construction due to Almgren and Pitts. Thus, if $L={\bf L}(\Pi)$, we deduce that ${\bf C}(S)$ must contain an element that is almost minimizing in annuli hence smooth. But note that ${\bf C}(S)$ could contain nonsmooth elements too. We are going to modify the combinatorial argument to obtain a critical sequence $S$ (meaning  ${\bf L}(S)={\bf L}(\Pi)$) such that every element of ${\bf C}(S)$ is almost smooth.

First we prove a lemma.
\subsection{Lemma}\label{lemma.convergence.thm}{\em Let $R, L>0$ and $m \in \mathbb{N}$. There exists $\eta=\eta(R,L,m)>0$ so that for every $V,T\in \mathcal V_n(M)$ with ${\bf M}(V)\leq 2L$, ${\bf M}(T)\leq 2L$, $V$ stationary,  $$V\llcorner  \left(M \setminus (\overline B_{2\eta}(p_1)\cup \cdots \cup \overline B_{2\eta}(p_t))\right)=T\llcorner  \left(M \setminus (\overline B_{2\eta}(p_1)\cup \cdots \cup \overline B_{2\eta}(p_t))\right)$$ for some collection $\{p_1,\dots,p_t\}\subset M$, $t\leq 3^{2m}$, and $$||V||(M)-\eta\leq ||T||(M)\leq||V||(M)+\eta,$$ then ${\bf F}(V,T)<R/2$. }
\begin{proof}
If false, we obtain  the existence of $V,T\in \mathcal V_n(M)$ and $\{p_1,\dots,p_{t'}\}\subset M$, $t'\leq 3^{2m}$, 
\begin{itemize}
\item[(i)] $||V||(M)=||T||(M)$ and $V=T$ on $M\setminus  (\overline B_{r}(p_1)\cup \cdots \cup \overline B_r(p_{t'}))$ for all $r>0$;
\item[(ii)]${\bf F}(V,T)\geq R/2$;
\item[(iii)] $V$ is stationary.
\end{itemize}
The Monotonicity formula for stationary varifolds implies the existence of a constant $C$ such that $||V||(B_r(q))\leq Cr^n$ for all $r$ sufficiently small. Thus we obtain from (i) that $V$ is identical to $T$, which contradicts (ii).
\end{proof}

The Isoperimetric Inequality of Federer and Fleming (Proposition 1.11 of \cite{almgren}), adapted to flat chains mod two,  implies there exists $\gamma_{iso}=\gamma(M)>0$ such that,  if $T \in \mathcal{Z}_n(M;\mathbb{Z}_2)$ satisfies ${\bf M}(T)\leq \gamma_{iso}$ then there exists $Q \in {\bf I}_{n+1}(M;\mathbb{Z}_2)$ with $\partial Q=T$ and ${\bf M}(Q)\leq {\bf M}(T)$.

Let $c=(3^m)^{3^m}$. The next theorem contains the essence of the combinatorial argument of \cite[Theorem 4.10]{pitts}.

\subsection{Theorem}\label{combinatorial.thm}{\em Consider a map $\phi: X(q)_0 \rightarrow \mathcal Z_n(M;\mathbb{Z}_2)$, and let $\mathcal{W}\subset \mathcal{V}_n(M)$ and $L = \sup_{x\in X(q)_0} {\bf M}(\phi(x))$.
Suppose  $R, \bar\varepsilon,\eta>0$ are such that $\bar\varepsilon<2R$ and if $x\in X(q)_0$ satisfies
$$
{\bf M}(\phi(x))\geq L-\bar\varepsilon {\rm \, \, and\, \,}  {\bf F}(|\phi(x)|, \mathcal W)\geq R,
$$
then there exist   $p(x)\in M$ and positive numbers $$r_1(x),\ldots, r_c(x), s_1(x),\ldots, s_c(x)$$ satisfying 
\begin{eqnarray*}
&&r_i(x)-2s_i(x)>2(r_{i+1}(x)+2s_{i+1}(x)), \, \, \, i=1,\dots,c-1,\\
&&r_1(x)+2s_1(x) < \eta, \\
&&r_c(x)-2s_c(x)>0,
\end{eqnarray*}
 such that $\phi(x)$ admits an $(\bar \varepsilon, \delta)$-deformation in each  $$A(p(x),r_j(x)-s_j(x),r_j(x)+s_j(x))\cap M,$$ $j=1,\ldots, c,$ for every $\delta>0$. If
 $$
 m{\bf f}(\phi)(1+4(3^m-1)s^{-1})<\min\{\frac{\bar\varepsilon}{3^{2m}8}, \gamma_{iso}\}
 $$
with  $s=\min_{x}\min\{s_1(x),\ldots, s_c(x)\}$,
then there exist $C=C(m,s)>0$,  an integer $k\geq q$,  and  a map $$\phi^*: X(k)_0 \rightarrow \mathcal Z_n(M;\mathbb{Z}_2),$$ such that  
\begin{itemize}
\item[(i)] $\phi^*$ is $(X,{\bf M})$-homotopic to $\phi$, through a discrete homotopy $$\psi: ([0,1] \times X)(k)_0 \rightarrow \mathcal{Z}_n(M;\mathbb{Z}_2)$$  with fineness less than or equal to $C{\bf f}(\phi),$
\end{itemize}
 and such that for every $t \in I(1,k)_0$, $x\in X(k)_0$, if $\hat{x} = {\bf n}(k+j,q+j)(x) \in X(q)_0$ then the following properties hold true:
 \begin{itemize}
\item[(ii)]  $$\phi^*(x)\llcorner (M\setminus \left(\overline B_{\eta}(p_1)\cup \cdots \cup \overline B_{\eta}(p_t)\right))=\phi(\hat x)\llcorner (M\setminus \left(\overline B_{\eta}(p_1)\cup \cdots \cup \overline B_{\eta}(p_t)\right))$$ for some collection $\{p_1,\dots,p_t\}\subset M$, $t\leq 3^m$;
\vskip 0.05in
\item[(iii)] ${\bf M}(\psi(t,x))\leq {\bf M}(\phi(\hat x))+2\cdot 3^{2m} m(1+4(3^m-1)s^{-1}){\bf f}(\phi)$;
\vskip 0.05in
\item[(iv)] $\phi^*(x)=\phi(\hat x)$ if ${\bf M}(\phi(\hat x))< L-\bar\varepsilon$ or if  $|\phi(\hat x)|\in {\bf B}^{\bf F}_{R}(\mathcal W)$;
\vskip 0.05in
\item[(v)] ${\bf M}(\phi^*(x)) \leq {\bf M}(\phi(\hat x))-\bar\varepsilon/4$ if  ${\bf M}(\phi(\hat x))\geq L-\bar\varepsilon/2$ and $|\phi(\hat x)|$ does not belong to ${\bf B}^{\bf F}_{2R}(\mathcal W);$
\item[(vi)] if ${\bf M}(\phi^*(x))\geq L-\bar\varepsilon/5$, then ${\bf F}(|\phi(\hat x)|, \mathcal W) \leq 2R.$ 
\end{itemize}
}

\begin{proof}

Recall that $X$ is a subcomplex of $I(m,j)$. We follow the arguments of Pitts  and explain the modifications one has to make to Parts 1 to 20 of \cite[Theorem 4.10]{pitts}.
\vskip0.05in

\noindent{{\bf Part 1}}:  Set $s=\min_{x}\min\{s_1(x),\ldots, s_c(x)\}$ and let
\begin{eqnarray*}
&&A_j(x)= A(p(x),r_j(x)-2s_j(x),r_j(x)+2s_j(x))\cap M,\\
&& a_j(x)= A(p(x),r_j(x)-s_j(x),r_j(x)+s_j(x))\cap M
 \end{eqnarray*}
 for all  $j=1,\ldots, c$, $x \in X(q)_0$ satisfying ${\bf M}(\phi(x))\geq L-\bar\varepsilon$ and  ${\bf F}(|\phi(x)|, \mathcal W)\geq R.$

\vskip0.05in

\noindent{{\bf Part 2}}:
 For each $\delta>0$, $j=1,\ldots, c$ and $x \in X(q)_0$ satisfying $$
{\bf M}(\phi(x))\geq L-\bar\varepsilon {\rm \, \, and\, \,}  {\bf F}(|\phi(x)|, \mathcal W)\geq R,
$$
 there exists a sequence $\phi(x)=T_1, T_2, \ldots, T_Q$ in $\mathcal Z_n(M;\mathbb{Z}_2)$ with ${\bf M}(T_Q)< {\bf M}(T_1)-\bar\varepsilon$ and such that
$$\s (T_i-T_1)\subset a_j(x),\quad {\bf M}(T_i-T_{i-1})\leq \delta,\quad {\rm \, and\, \,} {\bf M}(T_i)\leq {\bf M}(T_1)+\delta $$
for all $1\leq i\leq Q$. We can suppose $Q=3^{\bar{N}}$ for some $\bar{N} \in \mathbb{N}$  by repeating the last element of the sequence if necessary.
\vskip0.05in

\noindent{{\bf Parts 3-7}}: Do not do anything.

 \vskip0.05in

\noindent{{\bf Part 8}}: Let $$\delta=2\,3^{2m}m{\bf f}(\phi)(1+4(3^m-1)s^{-1}).$$  

\vskip0.05in

\noindent{{\bf Part 9}}: 
We  claim the existence of a  positive integer $N_1$ such that whenever $\sigma$ is a  cell of $X(q)$,  $\{x_1,\dots,x_{3^m}\}\subset \sigma_0$, $\mu\in\{1,\ldots,c\}$,  $$\min_{x\in\sigma_0}\{{\bf M}(\phi(x)\}\geq L-\bar\varepsilon, \quad\mbox{and}\quad \min_{x\in\sigma_0}\{{\bf F}(|\phi(x)|,\mathcal W)\}\geq R,$$ 
then for all  $j=1,\dots,3^m$, there is a sequence $\{T(j,l)\}_{l=1}^{3^{N_1}}$ in $\mathcal Z_n(M;\mathbb{Z}_2)$  such that 
\begin{itemize}
\item $T(j,1)=\phi(x_j)$;
\vskip0.02in
\item $\bigcup_{l=1}^{3^{N_1}}\s (T(j,l)-\phi(x_j))\subset A_\mu(x_1);$
\vskip0.02in

\item ${\bf M}(T(j,l)-T(j,l-1))<3^{-2m}\delta$ for all $l=2,\ldots, 3^{N_1}$;
\vskip0.02in

\item ${\bf M}(T(j,l)-T(j',l))<3^{-2m}\delta$ for all $l=1,\ldots, 3^{N_1}$ and $j'=1,\dots,3^m$;
\vskip0.02in

\item ${\bf M}(T(j,l))\leq {\bf M}(\phi(x_j))+3^{-2m}\delta$ all $l=1,\ldots, 3^{N_1}$;
\vskip0.02in

\item $||T(j,3^{N_1})||(A_\mu(x_1))< ||\phi(x_j)||(A_\mu(x_1))-\bar\varepsilon/2.$
\end{itemize}
\vskip0.05in

The construction of $T(j,l)$ follows  like in \cite[pages 167-168]{pitts}. The construction of $T(j,2)$ is so that
$$
T(j,2) \llcorner a_\mu(x_1)=\phi(x_1)  \llcorner a_\mu(x_1).
$$

There exist $Q_2,\dots,Q_{3^m}$ in ${\bf I}_{n+1}(M,\mathbb{Z}_2)$ such that
$$
\partial Q_j = \phi(x_j)-\phi(x_1),
$$
and
$$
{\bf M}(Q_j) \leq M(\phi(x_j)-\phi(x_1)) \leq m {\bf f}(\phi).
$$
If $d(y)=|y-p(x_1)|$, we can use Lemma 28.5 of \cite{simon} to  find $t_1,t_2 \in \mathbb{R}$ such that
\begin{eqnarray*}
&&r_\mu(x_1)+s_\mu(x_1) < t_1 < r_\mu(x_1)+2s_\mu(x_1), \\
&& r_\mu(x_1)-2s_\mu(x_1) < t_2 < r_\mu(x_1)-s_\mu(x_1),
\end{eqnarray*}
$$
||\phi(x_j)||(\partial B(p(x_1),t_1) \cup \partial B(p(x_1),t_2))=0, \, \, j=1, \dots, 3^m,
$$
$$
\langle Q_j,d,t_1 \rangle \in {\bf I}_n(M,\mathbb{Z}_2), \, \, \langle Q_j,d,t_2 \rangle \in {\bf I}_n(M,\mathbb{Z}_2), \, \, j=2, \dots, 3^m,
$$
$$
{\bf M}(\langle Q_j,d,t_1 \rangle) \leq (3^m-1) m {\bf f}(\phi)s^{-1}, \, \, j=2, \dots, 3^m,
$$
$$
{\bf M}(\langle Q_j,d,t_2 \rangle) \leq (3^m-1) m {\bf f}(\phi)s^{-1}, \, \, j=2, \dots, 3^m.
$$

For each $j=2,\dots,3^m$, we define
\begin{eqnarray*}
&&T(j,2)  = \phi(x_j) \llcorner (M \setminus B(p(x_1),t_1)) + \phi(x_1) \llcorner (M \cap A(p(x_1),t_2,t_1))\\
&&\hspace{2cm}  - \langle Q_j,d,t_1 \rangle + \langle Q_j,d,t_2 \rangle +  \phi(x_j) \llcorner (M \cap B(p(x_1),t_2)).
\end{eqnarray*}

Then
$$
T(j,2) \in \mathcal{Z}_n(M,\mathbb{Z}_2),
$$
$$
{\bf M}(T(j,2) -T(j,1))\leq m{\bf f}(\phi)(1+2(3^m-1)s^{-1}),
$$
$$
{\bf M}(T(j,2)-T(j',2)) \leq m{\bf f}(\phi)(1+4(3^m-1)s^{-1}), \, j'=1,\dots, 3^m,
$$
$$
m{\bf f}(\phi)(1+4(3^m-1)s^{-1}) < 3^{-2m}\delta.
$$

Now 
$$T(j,2)  \llcorner a_\mu(x_1) = \phi(x_1) \llcorner a_\mu(x_1)$$ for every $j=1,\dots, 3^m$.

By Part 2,  there exists a sequence $S_2=\phi(x_1), S_3, \dots, S_{3^{N_1}}$ in $\mathcal Z_n(M;\mathbb{Z}_2)$ with ${\bf M}(S_{3^{N_1}})< {\bf M}(\phi(x_1))-\bar \varepsilon$ and such that
$$\s (S_l-\phi(x_1))\subset a_\mu(x_1),$$
$$ {\bf M}(S_l-S_{l-1})< \frac12 3^{-2m}\delta, \, \, {\bf M}(S_l)< {\bf M}(\phi(x_1))+\frac12 3^{-2m}\delta$$
for all $3\leq l\leq 3^{N_1}$.

We set
$$
T(j,l) = T(j,2) \llcorner (M \setminus a_\mu(x_1))+ S_l \llcorner a_\mu(x_1)
$$
for $l=3,\dots, 3^{N_1}$.

\vskip0.05in

\noindent{{\bf Part 10}}: We choose a positive integer $N_2\geq q$ and define 
$$ f_2:X(N_2)_0\rightarrow \{0, 1,\ldots,3^{N_1}\}$$
such that, with $\hat\phi=\phi\circ{\bf n}(N_2+j,q+j)$, we have
\begin{itemize}
\item $f_2(x)=3^{N_1}$ whenever $x\in X(N_2)_0$, $${\bf M}(\hat\phi(x))\geq L-\bar\varepsilon/2,\quad\mbox{and}\quad{\bf F}(|\hat\phi(x)|,\mathcal W)\geq 2R;$$
\vskip 0.05in
\item if $x,y\in X(N_2)_0$ belong to a common cell and $f_2(x)>0$, then ${\bf M}(\hat\phi(y))\geq L-\bar\varepsilon$ and
${\bf F}(|\hat\phi(y)|,\mathcal W)\geq R;$
\vskip0.05in
\item $f_2(x)=0$ whenever $x\in X(N_2)_0$ and $${\bf M}(\hat\phi(x))< L-\bar\varepsilon\quad\mbox{or}\quad{\bf F}(|\hat\phi(x)|,\mathcal W)< R;$$
\vskip0.05in
\item $|f_2(x)-f_2(y)|\leq 1$ whenever $x,y\in X(N_2)_0$ belong to a common cell.
\end{itemize}
This can be seen in the following way. Let $A\subset X$, $B\subset X$ be, respectively, the union of all cells of $X(q)$ that contain a vertex $x$ such that 
$${\bf M}(\phi(x))\geq L-\bar\varepsilon/2\quad\mbox{and}\quad{\bf F}(|\phi(x)|,\mathcal W)\geq 2R,$$
or
$${\bf M}(\phi(x))< L-\bar\varepsilon\quad\mbox{or}\quad{\bf F}(|\phi(x)|,\mathcal W)<R.$$
Then  $A$ and $B$ are disjoint closed sets. Choose a continuous  function $h:X\rightarrow [0,3^{N_1}]$ such that $h$ is $3^{N_1}$ on $A$ and $0$ on $B$. Given $N_2\geq q$ and $x\in X(N_2)_0$ we set $f_2(x)$ to be the integer part of $h(x)$. For $N_2$ sufficiently large, $f_2$ satisfies the desired properties. 
\vskip0.05in

\noindent{{\bf Part 11}}: We choose
$$
f_3:X(N_2) \rightarrow X(N_2)_0
$$
so that $f_3(\sigma) \in \sigma_0$.

\vskip0.05in

\noindent{{\bf Part 12}}: We say $\sigma \in X(N_2)$ is a {\it bad cell} if $f_2(x)>0$ for some $x\in \sigma_0$, and is a {\it good cell} otherwise. We denote by $X_b(N_2), X_g(N_2)$ the set of bad and good cells, respectively. 


Note that if $\sigma$ is a bad cell then, by Part 10, 
$${\bf M}(\hat\phi(y))\geq L-\bar\varepsilon\quad\mbox{and}\quad{\bf F}(|\hat\phi(y)|,\mathcal W)\geq R$$
for every $y \in \sigma_0$.

For each cell $\sigma \in X_b(N_2)$, let $f_4(\sigma)={\bf n}(N_2+j,q+j)(f_3(\sigma))$. Then
$${\bf M}(\phi(f_4(\sigma)))\geq L-\bar\varepsilon\quad\mbox{and}\quad{\bf F}(|\phi(f_4(\sigma))|,\mathcal W)\geq R.$$
Let $a(\sigma)$ be the collection of disjoint annuli $A_j(f_4(\sigma))$, $j=1, \dots, c$.


\vskip0.05in

\noindent{{\bf Part 13}}: It follows from Proposition 4.9 of \cite{pitts} that we can choose
$$
f_5: X_b(N_2)  \rightarrow 2^M
$$
such that
\begin{itemize}
\item $f_5(\tau)$ is the annulus $A_j(f_4(\tau))$ for some $j \in \{1,2,\dots,c\}$, for every $\tau \in X_b(N_2)$,
\item the annuli $f_5(\tau)$ and $f_5(\sigma)$ are disjoint whenever $\sigma \neq \tau$ and $\sigma$ and $\tau$ are faces of a common cell in $I(m,N_2+j)$.
\end{itemize}
Moreover, $\# {\rm \, image\,}(f_5) \leq 3^{2m}.$

\vskip0.05in

\noindent{{\bf Part 14}}: Let $\tau \in X(N_2)$ and $\mu \in \{0,\dots,3^{N_1}\}$. For each $x\in \tau_0$, we define $f_6(\tau,x,\mu) \in \mathcal{Z}_n(M,\mathbb{Z}_2)$ as follows.

If $\tau$ is a good cell, we set $f_6(\tau,x,\mu) =\hat\phi(x)= \phi({\bf n}(N_2+j,q+j)(x))$.

If $\tau$ is a bad cell, then by Part 9, for each $y \in \tau_0$, there exists a sequence 
$$
T_y(1)= \phi({\bf n}(N_2+j,q+j)(y)), T_y(2), \dots, T_y(3^{N_1})
$$
in $\mathcal{Z}_n(M,\mathbb{Z}_2)$ such that
$$
{\rm spt\,}(T_y(1)-T_y(l)) \subset f_5(\tau), \, \, l=1,\dots, 3^{N_1},
$$
$$
{\bf M}(T_y(l)-T_y(l-1))<3^{-2m}\delta, \, \, l=2,\dots, 3^{N_1},
$$
$$
{\bf M}(T_y(l)-T_z(l)) < 3^{-2m}\delta, \, \, l=1, \dots, 3^{N_1}, z\in \tau_0,
$$
$$
{\bf M}(T_y(l)) \leq {\bf M}(T_y(1)) + 3^{-2m}\delta,  \, \, l=1, \dots, 3^{N_1},
$$
$$
{\bf M}(T_y(3^{N_1})) < M(T_y(1))-\bar\varepsilon/2.
$$

We set
$$
f_6(\tau,x,0) = \phi({\bf n}(N_2+j,q+j)(x)),
$$
$$
f_6(\tau,x,\mu) = T_x(\mu) \, {\rm if\,} 1\leq \mu \leq f_2(x),
$$
and
$$
f_6(\tau,x,\mu) = T_x(f_2(x)) \, {\rm if\,} f_2(x)\leq \mu\leq 3^{N_1}.
$$

\vskip0.05in

\noindent{{\bf Part 15}}: Let $N_3=N_1+N_2+2$, and ${\bf n}={\bf n}(N_3+j,N_2+j)$.

\vskip0.05in

\noindent{{\bf Part 16}}: We define
$$
f_7: X(N_3)_0 \rightarrow X(N_2)
$$
so that for each $x \in X(N_3)_0$, $f_7(x)$ is the unique cell of least dimension in $X(N_2)$ that contains $x$.

\vskip0.05in

\noindent{{\bf Part 17}}: We choose
$$
f_8: \{(x,\tau): x\in X(N_3)_0, \tau {\rm \, is \, a \, face \, of\,} f_7(x)\} \rightarrow \{0,1,\dots, 3^{N_1}\}
$$
so that it satisfies the following properties:
\begin{itemize}
\item for each $x \in X(N_3)_0$, there exists a face $\tau$ of $f_7(x)$ with ${\bf n}(x) \in \tau_0$ and $f_8(x,\tau)=3^{N_1}$,
\item $|f_8(x,\tau)-f_8(y,\tau)| \leq {\bf d}(x,y)$ for all $(x,\tau), (y,\tau) \in {\rm dmn\,}(f_8)$,
\item if $(x,\tau)\in {\rm dmn\,}(f_8)$, $y \in X(N_3)_0$ and $(y,\tau)\notin   {\rm dmn\,}(f_8)$, then $$f_8(x,\tau) \leq {\bf d}(x,y),$$
\item if $(x,\tau) \in {\rm dmn}(f_8)$ and $f_8(x,\tau)>0$, then  ${\bf n}(x) \in \tau_0$.
\end{itemize}

The construction of $f_8$ is just as in Part 17 of \cite{pitts}.





\vskip0.05in

\noindent{{\bf Part 18}}: We define
$$
\psi: I(1,N_3)_0 \times X(N_3)_0 \rightarrow \mathcal{Z}_n(M,\mathbb{Z}_2)
$$
as follows.

Let $x \in X(N_3)_0$. If $f_2({\bf n}(x))=0$, we set
$$
\psi(j,x) = \phi \circ {\bf n}(N_2+j,q+j) \circ {\bf n}(x)
$$
for each $j \in I(1, N_3)_0$.

Suppose $f_2({\bf n}(x))>0$. Let $F_x$ be the set of all  $\tau \in X(N_2)$ such that $\tau$  is a face of $f_7(x)$ with  ${\bf n}(x) \in \tau_0$. Every $\tau \in F_x$ is a bad cell, hence the annulus $f_5(\tau)$ is well-defined. Given $j \in I(1,N_3)_0$, we will define $\psi(j,x)$ by giving its restrictions to each such annulus and to the complement of their union
$$Z=M \setminus \cup \{f_5(\tau): \tau \in F_x\}.$$

We set
$$
\psi(j,x) \llcorner Z = \phi\circ {\bf n}(N_2+j,q+j) \circ {\bf n}(x) \llcorner Z
$$
for every $0 \leq j \leq 1$.

Then we put
$$
\psi(j,x)  \llcorner f_5(\tau) = \phi \circ {\bf n}(N_2+j,q+j) \circ {\bf n}(x) \llcorner f_5(\tau)
$$
if $j =0$ or $3^{-N_3}$;
$$
\psi(j,x)  \llcorner f_5(\tau)=f_6(\tau, {\bf n}(x), 3^{N_3} \cdot j) \llcorner f_5(\tau)
$$
if $1 \leq j\cdot 3^{N_3} \leq \min\{f_2({\bf n}(x)), f_8(x,\tau)\},$
and
$$
\psi(j,x)  \llcorner f_5(\tau)=f_6(\tau, {\bf n}(x), \min\{f_2({\bf n}(x)), f_8(x,\tau)\}) \llcorner f_5(\tau)
$$
if $\min \{f_2({\bf n}(x)), f_8(x,\tau)\} \leq j \cdot 3^{N_3} \leq 3^{N_3}$.


We would like to estimate the fineness of $\psi$ in terms of the fineness of $\phi$.

We have, by Part 14, that
\begin{eqnarray*}
{\bf M}(\psi(j,x)-\psi(j-1,x)) < \sum_{f_5(\tau), \tau \in F_x} 3^{-2m} \delta \leq \delta.
\end{eqnarray*}

Let $x,x' \in X(N_3)_0$ such that ${\bf d}(x,x')=1$. Then the 1-cell $[x,x'] \in X(N_3)$ and there is a top  cell $\sigma \in X(N_2)$ such that
$[x,x'] \subset \sigma$. Then both $f_7(x)$ and $f_7(x')$ are faces of $\sigma$. Notice that also ${\bf d}({\bf n}(x),{\bf n}(x'))\leq 1$, in particular $|f_2({\bf n}(x))-f_2({\bf n}(x'))| \leq 1$.

If $f_2({\bf n}(x))=0$, then $f_2({\bf n}(x')) \leq 1$ and
\begin{eqnarray*}
&&{\bf M}(\psi(j,x)-\psi(j,x'))\\
&&={\bf M}(\phi \circ {\bf n}(N_2+j,q+j) \circ {\bf n}(x)-\phi \circ {\bf n}(N_2+j,q+j) \circ {\bf n}(x'))
 \leq {\bf f}(\phi).
\end{eqnarray*}
The same inequality holds if $f_2({\bf n}(x'))=0$.

Hence we can suppose both $f_2({\bf n}(x))$ and $f_2({\bf n}(x'))$ are positive. Suppose  $\tau \in F_x \setminus  F_{x'}$. Then $\tau$ is a face of $f_7(x)$ and ${\bf n}(x)\in \tau_0$. If $\tau$ is not a face of $f_7(x')$, then $(x',\tau) \notin {\rm dmn}(f_8)$. Since $(x,\tau) \in {\rm dmn}(f_8)$, by Part 17 we have that $f_8(x,\tau) \leq {\bf d}(x,x')=1$. If $\tau$ is a face of $f_7(x')$ but ${\bf n}(x') \notin \tau_0$, Part 17 implies
$f_8(x',\tau)=0$. We conclude, using Part 17,  that if $\tau \in F_x \setminus  F_{x'}$ then
$$
f_8(x,\tau) \leq 1.
$$
Similarly, if $\tau \in F_{x'} \setminus  F_{x}$ then
$
f_8(x',\tau) \leq 1.
$

Hence
\begin{eqnarray*}
&&{\bf M}(\psi(j,x)-\psi(j,x')) \leq \sum_{\tau \in F_x\setminus F_{x'}} {\bf M}((\psi(j,x)-\psi(j,x')) \llcorner f_5(\tau)) \\
&&+  \sum_{\tau \in F_{x'}\setminus F_{x}} {\bf M}((\psi(j,x)-\psi(j,x')) \llcorner f_5(\tau))\\
&& +  \sum_{\tau \in F_{x}\cap F_{x'}} {\bf M}((\psi(j,x)-\psi(j,x')) \llcorner f_5(\tau))\\
&& + {\bf M}\left((\psi(j,x)-\psi(j,x')) \llcorner \left(M \setminus \cup_{\tau \in F_x \cup F_{x'}} f_5(\tau)\right)\right)\\
&&\leq {\bf f}(\phi)+2\delta.\\
\end{eqnarray*}

Hence
$$
{\bf f}(\psi) \leq {\bf f}(\phi)+2\delta=C(m,s)\, {\bf f}(\phi).
$$

\vskip0.05in

\noindent{{\bf Part 19}}:  Define $$\phi^*: X(N_3)_0 \rightarrow \mathcal{Z}_n(M,\mathbb{Z}_2)$$ by
$$
\phi^*(x)=\psi(1,x).
$$

\vskip0.05in

\noindent{{\bf Part 20}}:  The fact that $\phi^*$ is $(X,{\bf M})$-homotopic to $\phi$ with fineness less than or equal to $C(m,s){\bf f}(\phi)$ follows from the fact that ${\bf f}(\psi)\leq C(m,s)\, {\bf f}(\phi)$ (Part 18).

Given $x \in X(N_3)_0$, let $\hat{x} = {\bf n}(N_2+j,q+j) \circ {\bf n}(x) \in X(q)_0$. Property (ii) follows from the fact that every annulus $f_5(\tau)$, $\tau \in F_x$, is contained in a ball of radius $\eta$ and the number of elements of $F_x$ is at most $3^{2m}$.  From Part 14 we get that ${\bf M}(\psi(t,x)) \leq {\bf M}(\phi(\hat{x})) + \delta$. Hence Property (iii) follows.

If ${\bf M}(\phi(\hat x))< L-\bar\varepsilon$ or if  $|\phi(\hat x)|\in {\bf B}^{\bf F}_{R}(\mathcal W)$, by Part 10 we have $f_2({\bf n}(x))=0$. Then, by Part 18, $\phi^*(x)=\phi(\hat{x})$. This implies Property (iv).

If ${\bf M}(\phi(\hat x))\geq L-\bar\varepsilon/2$ and   ${\bf F}(|\phi(\hat x)|,  \mathcal{W})\geq 2R$, we have by Part 10 that  $f_2({\bf n}(x))=3^{N_1}$. By Part 17, there is  a face $\tilde{\tau}\in F_x$  with $f_8(x,\tilde{\tau})=3^{N_1}$. Hence
\begin{eqnarray*}
&&{\bf M}(\phi^*(x)) = \sum_{\tau \in F_x} {\bf M}(\phi^*(x) \llcorner f_5(\tau)) + {\bf M}(\phi^*(x) \llcorner \left(M \setminus \cup_{\tau \in F_x} f_5(\tau)\right))\\
&&\leq {\bf M}(\phi^*(x) \llcorner f_5(\tilde{\tau})) +\sum_{\tau \in F_x, \tau \neq \tilde{\tau}} {\bf M}(\phi^*(x) \llcorner f_5(\tau)) \\
&&+ {\bf M}(\phi(\hat{x}) \llcorner \left(M \setminus \cup_{\tau \in F_x} f_5(\tau)\right))\\
&&\leq{\bf M}(\phi(\hat{x}) \llcorner f_5(\tilde{\tau}))-\bar\varepsilon/2 + \sum_{\tau \in F_x, \tau \neq \tilde{\tau}} {\bf M}(\phi(\hat{x}) \llcorner f_5(\tau)) + \sum_{\tau \in F_x, \tau \neq \tilde{\tau}} 3^{-2m}\delta\\
&&+ {\bf M}(\phi(\hat{x}) \llcorner \left(M \setminus \cup_{\tau \in F_x} f_5(\tau)\right))\\
&&\leq{\bf M}(\phi(\hat{x}))-\bar\varepsilon/2+ \delta.
\end{eqnarray*}
This proves Property (v).

Property (vi) follows from Properties (iii) and (v).

\end{proof}

Suppose again that $X$ is a $k$-dimensional cubical subcomplex of $I(m,j)$. Let  $\Phi:X \rightarrow \mathcal Z_n(M^{n+1};{\bf F};\Z_2)$ be a continuous map. We let $\Pi$ be the class of all continuous maps $$\Phi':X \rightarrow \mathcal Z_n(M^{n+1};{\bf F};\Z_2)$$ such that $\Phi$ and $\Phi'$ are homotopic to each other in the flat topology. 
Let $L={\bf L}(\Pi)$ and $\mathcal W_L$ be  the set of all  stationary integral varifolds  in $M$ with mass equal to $L$  and whose support is a closed, minimal hypersurface that is smooth and embedded outside a set of Hausdorff dimension $n-7$.

The previous theorem has the following continuous counterpart.
\subsection{Theorem}\label{combinatorial.thm.cont}{\em Consider a pulled-tight minimizing sequence $\{\Phi_i\}\in\Pi$ of continuous maps in the mass topology
$$\Phi_i: X \rightarrow \mathcal Z_n(M;{\bf M};\mathbb{Z}_2),  \quad i\in\N,$$
with $L={\bf L}(\Pi)={\bf L}(\{\Phi_i\})$.

For all $R,\delta, \bar\delta>0$, there exist  $0<\eta<R$, $0<\delta_1<\bar \delta$, $\delta_2>0$, $N\in\N$,  $C>0$ and a sequence of continuous maps in the mass topology $$\Psi_i: X \rightarrow \mathcal Z_n(M;{\bf M};\mathbb{Z}_2), \quad i\geq N$$ such that for all $i\geq N$
\begin{itemize}
\item[(i)] $\Psi_i$ is homotopic to $\Phi_i$ in the flat topology;
\item[(ii)]  for all $x\in X$ there exist $T\in \mathcal Z_n(M;\mathbb{Z}_2)$ and a set  $\{p_1,\dots,p_t\} \subset  M$, $t \leq 3^{2m}$, such that
$${\bf M}( (\Phi_i(x)-T) \llcorner (M\setminus  (\overline B_{\eta}(p_1)\cup \cdots \cup \overline B_{\eta}(p_t))   ))\leq 2^{-i}\quad\mbox{and}\quad {\bf F}(|T|,|\Psi_i(x)|)\leq \delta;$$
\vskip 0.05in
\item[(iii)] ${\bf M}(\Psi_i(x))\leq {\bf M}(\Phi_i(x))+C\cdot 4^{-i}$ for all $x\in X$;
\vskip 0.05in
\item[(iv)] ${\bf F}(|\Psi_i(x)|,|\Phi_i(x)|)\leq \delta$ if ${\bf M}(\Phi_i(x))< L-\delta_1$;
\vskip 0.05in
\item[(v)] $\{\Psi_i\}$ is a pulled-tight minimizing sequence and  
$$ {\bf M}(\Psi_i(x))\geq L-\delta_2 \implies |\Psi_i(x)| \in {\bf B}^{\bf F}_{3R}(\mathcal W_L\cap{\bf C}(\{\Phi_i\})).$$
\end{itemize}
}

\begin{proof}
 Choose   a nondecreasing sequence $\{n_i\}\subset \mathbb{N}$ such that
$$\sup\{{\bf M}(\Phi_i(x)-\Phi_i(y)): x,y \mbox{ in a common cell of }X(n_i)\}\leq 4^{-i}.$$
Consider the sequence $S=\{\phi_i\}$, $\phi_i:X(n_i)_0 \rightarrow  Z_n(M;\mathbb{Z}_2)$, with   $\phi_i(x)=\Phi_i(x)$ for every $x \in X(n_i)_0$. It satisfies ${\bf f}(\phi_i)\leq 4^{-i}$, $L={\bf L}(S)$, and ${\bf C}(S)={\bf C}(\{\Phi_i\})$. If $i$ is sufficiently large, $\phi_i$ admits an Almgren extension $\tilde{\Phi}_i:X \rightarrow \mathcal{Z}_n(M;{\bf M};\mathbb{Z}_2)$ (by Theorem 3.10 of \cite{marques-neves-infinitely}) that is homotopic to $\Phi_i$ in the mass topology by Proposition \ref{close.implies.homotopic}.

Let  $\eta=\eta(R,L,m)>0$ as in Lemma \ref{lemma.convergence.thm}.
An immediate consequence of Theorem \ref{regularity.thm} is that for all  $V\in {\bf C}(S)\setminus \mathcal W_L$ there exists $p\in M$ so that for all $r<\eta/2$ there is $s<r$ such that $V$ is not $\Z_2$-almost minimizing in $A(p,s,r)\cap M$. 

Set $\mathcal W = \mathcal W_L \cap {\bf C}(S)$ and $\mathcal K={\bf C}(S)\setminus {\bf B}^{\bf F}_{R}(\mathcal W).$ Therefore for all $V\in \mathcal K$ we obtain the existence of $p\in \s\, V$ and positive numbers $r_1,\ldots, r_c$, $s_1,\ldots, s_c$ satisfying 
\begin{eqnarray*}
&&r_i-2s_i>2(r_{i+1}+2s_{i+1}), \, \, \, i=1,\dots,c-1,\\
&&r_c-2s_c>0,
\end{eqnarray*}
 such that $V$ is not $\Z_2$-almost minimizing in $A(p,r_j-s_j,r_j+s_j)\cap M$ and $A(p,r_j-2s_j,r_i+2s_j)\subset B_\eta(p)$, for all $j=1,\ldots, c.$

Set $s(V)=\min\{s_1,\ldots, s_c\}$ and let
$$A_j(V)= A(p,r_j-2s_j,r_j+2s_j)\cap M,\quad a_j(V)= A(p,r_j-s_j,r_j+s_j)\cap M$$
 for all  $j=1,\ldots, c$.
 
  From the fact that, for each $j=1,\ldots, c$, $V\in \mathcal K$ is not $\Z_2$-almost minimizing in $a_j(V)$,
  there is $\varepsilon(V)>0$ with the property that for each $\delta>0$ and each $T\in\mathcal Z_n(M;\mathbb{Z}_2)$ with ${\bf F}(V,|T|)<\varepsilon(V)$, there exists a sequence $T=T_1, T_2, \ldots, T_q$ in $\mathcal Z_n(M;\mathbb{Z}_2)$ with ${\bf M}(T_q)< {\bf M}(T)-\varepsilon(V)$ and such that
$$\s (T_i-T)\subset a_j(V),\quad {\bf M}(T_i-T_{i-1})\leq \delta,\quad {\rm \, and\, \,} {\bf M}(T_i)\leq {\bf M}(T)+\delta $$
for all $1\leq i\leq q$. 

Because $\mathcal K$ is compact there is a finite set $\{V_i\}_{i=1}^{\nu}$ in $\mathcal{K}$ such that
$$\mathcal K\subset \bigcup_{i=1}^{\nu}{\bf B}^{\bf F}_{\frac{\varepsilon(V_i)}{4}}(V_i).$$

We set $\varepsilon_1=\min\{\varepsilon(V_i)\}_{i=1}^{\nu}$,
\begin{multline*}
\bar\varepsilon=\min\left\{\frac R 3 , \frac{\varepsilon_1}{2}, \frac{\bar\delta}{4}, \frac 1 4 \sup\{\varepsilon:{\bf F}(V,V_i)<\varepsilon(V_i)/4 \mbox{ for some }i=1,\ldots,\nu,\right.\\
\mbox{whenever }V\in {\bf K}(S)\setminus {\bf B}^{\bf F}_{R}(\mathcal W)\mbox{ and }{\bf M}(V)\geq L-2\varepsilon\}\bigg\},
\end{multline*}
and $s=\min\{s(V_j):j=1,\ldots,c\}$. 

Choose a positive integer $N$ so that whenever $i\geq N$ the following four properties hold:
 \begin{itemize}
 \item If $x\in X(n_i)_0$ is such that $${\bf M}(\phi_i(x))\geq \sup_{x\in X(n_i)_0} {\bf M}(\phi_i(x)) -\bar\varepsilon\quad\mbox{and}\quad{\bf F}(|\phi_i(x)|,\mathcal W)\geq R,$$ then ${\bf F}(|\phi_i(x)|, V_j)<\varepsilon(V_j)/2$ for some $j=1,\ldots,\nu$;
 \item $m{\bf f}(\phi_i)(1+4(3^m-1)s^{-1})<\min\{\frac{\bar\varepsilon}{3^m8},\gamma_{iso}\}.$
 \end{itemize}

From Theorem \ref{combinatorial.thm}, we obtain a sequence $\{\phi^*_i\}_{i\geq \bar N}$, $$\phi_i^*:X(k_i)_0 \rightarrow  Z_n(M;\mathbb{Z}_2),$$  with fineness less than or equal to $C(m,s){\bf f}(\phi_i)$ that satisfies Theorem  \ref{combinatorial.thm} (i)--(vi).

Choose $N_1\geq N$   so that Theorem 3.10 of \cite{marques-neves-infinitely} can be applied to all $\phi^*_i$ with $i\geq N_1$.  We denote by $ \Phi^*_i$ the Almgren extension of $\phi^*_i$ that satisfies, by Theorem 3.10 of  \cite{marques-neves-infinitely},
\begin{eqnarray*}
&&\sup\{{\bf M}(\Phi^*_i(x)-\Phi^*_i(y)): x,y\mbox{ in a common cell of }\dmn(\phi^*_i )\}\\
&&\leq C_0C(m,s){\bf f}(\phi_i)
\end{eqnarray*}
for all $i\geq  N_1$. Then $\Phi_i^*$ is homotopic to $\tilde{\Phi}_i$ in the flat topology for sufficiently large $i$, by Proposition 3.11 of \cite{marques-neves-infinitely}.

Given $x' \in X$, there exists $x \in X(k_i)_0$ such that $x, x'$ belong to a common cell of $X(k_i)$ and ${\bf M}(\Phi_i^*(x') - \phi_i^*(x)) \leq C_0 C(m,s){\bf f}(\phi_i)$. Let $T = \phi_i^*(x)$. Let $\hat{x}={\bf n}(k_i+j,n_i+j)(x) \in X(n_i)_0$ be as in Theorem \ref{combinatorial.thm} so by part (ii) of that theorem we have
$$
T\llcorner (M\setminus  (\overline B_{\eta}(p_1)\cup \cdots \cup \overline B_{\eta}(p_t))   ) =  \phi_i(\hat{x}) \llcorner (M\setminus  (\overline B_{\eta}(p_1)\cup \cdots \cup \overline B_{\eta}(p_t))   )
$$
for some $\{p_1,\dots,p_t \} \subset M$, $t \leq 3^{2m}$. Since $x'$ and $\hat{x}$ belong to a common cell in $X(n_i)$, we have that ${\bf M}(\Phi_i(x') - \Phi_i(\hat{x}))\leq 4^{-i}$. Notice that $\Phi_i(\hat{x}) = \phi_i(\hat{x})$. Hence
$${\bf M}( (\Phi_i(x')-T) \llcorner (M\setminus  (\overline B_{\eta}(p_1)\cup \cdots \cup \overline B_{\eta}(p_t))   ))\leq 4^{-i}.$$
From the above, we also have ${\bf M}(\Phi_i^*(x') - T) \leq C_0C(m,s) {\bf f}(\phi_i)$. 

We have
\begin{eqnarray}\label{massbound.combinatorial.thm.cont}
{\bf M}(\Phi_i^*(x')) &\leq& {\bf M}(\phi_i^*(x)) + C_0C(m,s){\bf f}(\phi_i) \nonumber \\
&\leq& {\bf M}(\phi_i(\hat{x}))+ C'(m,s) {\bf f}(\phi_i)\nonumber \\
&\leq& {\bf M}(\Phi_i(x')) + 4^{-i} + C'(m,s) {\bf f}(\phi_i).
\end{eqnarray}
Also 
$${\bf M}((\Phi_i(x')-\Phi_i^*(x')) \llcorner (M \setminus (\overline B_{\eta}(p_1)\cup \cdots \cup \overline B_{\eta}(p_t))))\leq 4^{-i} +C_0C(m,s){\bf f}(\phi_i).$$

Choose $N_2\geq N_1$ so that $4^{-N_2}\leq \bar\varepsilon/2$. We  deduce from Theorem  \ref{combinatorial.thm} (iv)  that
$${\bf M}(\Phi_i(x'))< L-2\bar\varepsilon \implies {\bf M}(\Phi_i(x')-\Phi^*_i(x'))\leq 4^{-i} + C_0C(m,s) {\bf f}(\phi_i),$$
for $i\geq N_2.$
From \eqref{massbound.combinatorial.thm.cont} we have that ${\bf L}(\{\Phi^*_i\}_{i\geq N_2})=L$. 

Let $W \in {\bf C}(\{\Phi^*_i\})$. Then $W$ is the varifold limit of $\phi_{i_l}^*(x_l)$ as $l \rightarrow \infty$, where $\{i_l\}\subset \{i\}$ and $x_l\in X(k_{i_l})_0$. Let $V$ be the varifold limit of $\phi_{i_l}(\hat x_l)$, after passing to a subsequence, where $\hat{x}_l={\bf n}(k_{i_l}+j,n_{i_l}+j)(x_l) \in X(n_{i_l})_0$. Property (iii) of Theorem \ref{combinatorial.thm} implies $||W||(M) \leq ||V||(M)$. But $||V||(M)\leq L$ and $L=||W||(M)$, hence $||V||(M)=||W||(M)=L$. Hence $V \in {\bf C}(\{\Phi_i\})$, and therefore $V$ is stationary. Since
$$
W\llcorner (M\setminus  (\overline B_{2\eta}(p_1)\cup \cdots \cup \overline B_{2\eta}(p_t))   ) =  V \llcorner (M\setminus  (\overline B_{2\eta}(p_1)\cup \cdots \cup \overline B_{2\eta}(p_t))   )
$$
for some $\{p_1,\dots,p_t \} \subset M$, $t \leq 3^{2m}$, Lemma \ref{lemma.convergence.thm} gives that ${\bf F}(V,W)<R/2$.
Theorem \ref{combinatorial.thm} (vi) implies ${\bf F}(V, \mathcal W)\leq 2R$.
Hence $${\bf C}(\{\Phi^*_i\})\subset \overline{\bf B}^{\bf F}_{5R/2}(\mathcal W=\mathcal W_L\cap{\bf C}(\{\Phi_i\})).$$ 

The minimizing sequence $\{\Phi_i^*\}$ is not necessarily pulled-tight. By Subsection \ref{pulltight}, with 
$\varepsilon=\delta/2$, we can find $\Psi^*_i$ continuous in the ${\bf F}$-topology such that $\{\Psi_i^*\}$ is a pulled-tight minimizing sequence for $\Pi$, ${\bf C}(\{\Psi_i^*\}) \subset {\bf C}(\{\Phi_i^*\})$, and  ${\bf M}(\Psi_i^*(x)) \leq {\bf M}(\Phi_i^*(x))$,  ${\bf F}(|\Psi_i^*(x)|, |\Phi_i^*|(x)) \leq \delta/2$ for every $x\in X$. By Proposition \ref{M.regularization}, we can find $\Psi_i$ continuous in the mass topology, homotopic to $\Psi_i^*$ in the ${\bf F}$-metric, with ${\bf F}(\Psi^*_i(x),\Psi_i(x))\leq 4^{-i}$ for all $i\geq N_2$, $x\in X$. In particular, $\Psi_i$ is homotopic to $\Phi_i$ in the flat topology for sufficiently large $i$. This proves Property (i).

Fix $i\geq N_2$, $x'\in X$. We can choose $T=\phi_i^*(x)$ as above so that
$${\bf M}( (\Phi_i(x')-T) \llcorner (M\setminus  (\overline B_{\eta}(p_1)\cup \cdots \cup \overline B_{\eta}(p_t))   ))\leq 4^{-i}.$$ Now
$${\bf F}(|\Psi_i(x')|,|T|) \leq {\bf F}(|\Phi_i^*(x')|,|T|) + 4^{-i} + \delta/2\leq C_0C(m,s){\bf f}(\phi_i) + 4^{-i} + \delta/2.$$ This proves Property (ii) of the theorem for large $i$.

We also have
\begin{eqnarray*}
{\bf M}(\Psi_i(x')) &\leq& {\bf M}(\Psi_i^*(x')) + 4^{-i}\\
&\leq& {\bf M}(\Phi_i^*(x')) + 4^{-i} \\
&\leq& {\bf M}(\Phi_i(x')) + 2 \cdot 4^{-i} + C'(m,s){\bf f}(\phi_i),
\end{eqnarray*}
hence Property (iii) of the theorem is proved.

If ${\bf M}(\Phi_i(x')) < L- 2\bar\varepsilon$, we have $ {\bf M}(\Phi_i(x')-\Phi^*_i(x'))\leq 4^{-i} + C_0 C(m,s){\bf f}(\phi_i)$. Therefore, in that case,
\begin{eqnarray*}
{\bf F}(|\Psi_i(x')|, |\Phi_i(x')|) &\leq& {\bf F}(|\Psi_i^*(x')|, |\Phi_i(x')|) +4^{-i}\\
&\leq&  {\bf F}(|\Phi_i^*(x')|, |\Phi_i(x')|) +4^{-i}+\delta/2\\
&\leq& 2\cdot 4^{-i}+\delta/2+  C_0 C(m,s){\bf f}(\phi_i).
\end{eqnarray*}
This proves Property (iv) for sufficiently large $i$.

Finally, $${\bf C}(\{\Psi_i\})={\bf C}(\{\Psi^*_i\})$$ and so $\{\Psi_i\}$ is a minimizing pulled-tight sequence. We also have
$$ {\bf C}(\{\Psi_i\})={\bf C}(\{\Psi^*_i\})\subset {\bf C}(\{\Phi^*_i\})\subset \overline {\bf B}^{\bf F}_{\frac{5R}{2}}(\mathcal W_L\cap{\bf C}(\{\Phi_i\})),$$ which implies  that Property  (v)  follows for some choice of $\delta_2$ and all $i$ sufficiently large.

\end{proof}

The next Corollary together with Theorem \ref{regularity.thm} implies Min-max Theorem \ref{minmax.continuous.thm}.

\subsection{Corollary}\label{combinatorial.cor}{\em Consider a  minimizing sequence $\{\Phi_i\}\subset\Pi$ of continuous maps in the ${\bf F}$-metric
$$\Phi_i: X \rightarrow \mathcal Z_n(M;{\bf F};\mathbb{Z}_2),  \quad i\in\N,$$
with $L={\bf L}(\Pi)={\bf L}(\{\Phi_i\})$. Then
$$
\mathcal W_L\cap{\bf C}(\{\Phi_i\})\neq \emptyset.
$$
}

\begin{proof}
By Section \ref{pulltight}, we can assume $\{\Phi_i\}$ is pulled-tight. Proposition \ref{M.regularization} implies that we can suppose $\Phi_i$ is continuous in the mass topology. If $\mathcal W_L\cap{\bf C}(\{\Phi_i\})= \emptyset$, Theorem \ref{combinatorial.thm.cont} gives a sequence $\{\Psi_i\} \subset \Pi$ that by part (v) of the theorem satisfies ${\bf M}(\Psi_i(x))<L-\delta_2$ for every $x\in X$ and $i$ sufficiently large. This is a contradiction, since $L={\bf L}(\Pi)$, and this proves the corollary.
\end{proof}

Now suppose $3 \leq (n+1) \leq 7$, so that the support of every element of $\mathcal W_L$ is smooth. We denote by $\mathcal W_{L,j}$, $\mathcal W_{L}^j$ the elements in $\mathcal W_L$ whose support has Morse index less than or equal to $j$  and bigger than or equal to $j$, respectively. A Riemannian metric $g$ is  called {\it bumpy} if there is no closed, smooth, immersed, minimal hypersurface that admits a non-trivial Jacobi  field. White showed in \cite{white2, white3} that bumpy metrics are generic in the usual $C^\infty$ Baire sense. Sharp's  Compactness Theorem \cite{sharp} implies $\mathcal W_{L,j}$ is a finite set if the metric $(M,g)$ is bumpy.

The next theorem proves the existence of minimizing sequences such that every element of the critical set is almost smooth.

\subsection{Theorem}\label{theorem.improv.minmax}{\em Suppose $(M,g)$ is a bumpy metric.  Consider  a minimizing sequence $\{\Psi_i\}\in \Pi$  such that $${\bf C}(\{\Psi_i\})\cap \mathcal W_{L}^{k+1}=\emptyset,$$ and let  $\Lambda_1={\bf C}(\{\Psi_i\})\cap  \mathcal W_{L}\subset  \mathcal W_{L, k}.$ For every $s>0$, there is a  minimizing sequence $\{\Phi_i\}\in \Pi$ such that ${\bf C}(\{\Phi_i\})\cap \mathcal{W}_L^{k+1}=\emptyset$ and for some $\gamma\in (0,L)$ and every sufficiently large $i$, we have
 $$\{|\Phi_i(x)|\in\mathcal V_n(M):{\bf M}(\Phi_i(x))\geq L-\gamma\}\subset \bigcup_{\Sigma\in \Lambda_1}{\bf B}_s^{\bf F}(\Sigma).$$
}
\medskip

\begin{proof}
Consider the finite set $\mathcal T$ of all elements in $\mathcal Z_n(M;\Z_2)$ whose support is contained in the support of some element in $\mathcal W_{L,k}$. We choose $s_1>0$ so that  for all $\Sigma\in \mathcal W_{L,k}$, $\Sigma'\in \mathcal T$,
$${\bf B}_{2s_1}^{\bf F}(\Sigma)\cap \mathcal W_{L,k} =\{\Sigma\}, \quad {\bf B}_{2s_1}^{\mathcal F}(\Sigma')\cap \mathcal T =\{\Sigma'\}.$$
There exists $0<s_2<s_1$ so  that for all $T\in \mathcal Z_n(M;\Z_2)$ with ${\bf F}(|T|,\Sigma)<2s_2$ for some $\Sigma\in \mathcal W_{L,k}$ there is $\Sigma'\in \mathcal T$ so that ${\mathcal F}(T,\Sigma')<s_1$. This is a consequence of the Constancy Theorem of flat chains and lower semicontinuity of mass with respect to flat topology. 

We can assume that $\{\Psi_i\}$ is pulled tight. By Proposition \ref{M.regularization}, we can also assume that each $\Psi_i$ is continuous in the mass topology. The set $\Lambda_1={\bf C}(\{\Psi_i\})\cap  \mathcal W_{L}\subset  \mathcal W_{L, k}$ is finite.

Applying Theorem \ref{combinatorial.thm.cont} to $\{\Psi_i\}\in \Pi$ with $R=s/16$, $\delta=\overline \delta=1$, we obtain  a pulled-tight minimizing sequence $\{\Psi^1_i\}\subset \Pi$ so that ${\bf C}(\{\Psi^1_i\})\subset {\bf B}_{s/4}^{\bf F}(\Lambda_1)$.  Using a slight variant of  the construction in Theorem 6.1 of  \cite{marques-neves-index}   (with $r=s/100$ in its notation), we obtain a pulled-tight min-max sequence $\{\Psi^2_i\}$ with ${\bf C}(\{\Psi^2_i\})\cap \mathcal W_{L}^{k+1}=\emptyset$ and  
\begin{equation}\label{critical.distance}
{\bf C}(\{\Psi^2_i\})\subset {\bf B}_{s/100}^{\bf F}({\bf C}(\{\Psi^1_i\}))\subset {\bf B}_{s/3}^{\bf F}(\Lambda_1).
\end{equation} We can also assume $\Psi_i^2$ is continuous in the mass topology for every $i$.

\end{proof}

The next proposition concerns sequences of mod 2 flat cycles such that the flat limit is not compatible to the varifold limit. Examples of such sequences are given in Section 2.4 of White \cite{white-currents}.

\subsection{Proposition}\label{flat.wrong} {\it Suppose $3\leq (n+1) \leq 7$. Let $$V = m_1 \cdot |\Sigma_1| + \cdots + m_N \cdot |\Sigma_N|$$ be an embedded minimal cycle, where $\{\Sigma_1, \dots, \Sigma_N\}$ is a disjoint collection, and 
$$\Sigma = n_1\cdot  \Sigma_1 + \cdots + n_N \cdot \Sigma_N,$$ $n_i \in \{0,1\}$, be a mod two flat cycle with 
$n_j \neq m_j {\rm \, mod\,} 2$ for some $1\leq j\leq N$.

For every $\alpha>0$,
there exists $\varepsilon>0$ such that if $T\in \mathcal{Z}_n(M;\mathbb{Z}_2)$ satisfies $${\bf F}(|T|, V) \leq \varepsilon {\rm \, \, and\, \,} \mathcal{F}(T, \Sigma)\leq \varepsilon,$$ then for every $p\in \Sigma_j$, $T$ admits an $(\varepsilon, \delta)$-deformation in 
$B(p,\alpha) \cap M$ for every $\delta>0$.
}

\begin{proof}
Let $r_0>0$ be sufficiently small so that, for every $p \in \s(V)$,
\begin{itemize}
\item the restriction of $v \mapsto |v-p|^2$ to $M \cap B(p,2r_0)$ is strictly convex, 
\item $\s(V)$ is transversal to $M \cap \partial B(p,s)$ for every $s< 2r_0$,
\item if $p\in \Sigma_l$, then $B(p,2r_0) \cap \Sigma_{l'}=\emptyset$ for every $1\leq l'\leq N$, $l'\neq l$.
\end{itemize}

Suppose the statement of the proposition is false, by contradiction.
This means there exist $T_i\in \mathcal{Z}_n(M;\mathbb{Z}_2)$ with ${\bf F}(|T_i|, V) \leq 1/i$ and $\mathcal{F}(T_i, \Sigma)\leq 1/i$ and  $p_i\in \Sigma_j$ such that
$T_i \in \mathfrak{a}(B(p_i,\alpha) \cap M;1/i,\delta_i)$ for some $\delta_i>0$. We can assume that $\alpha\leq r_0$, that $\{p_i\}$ converges to $p \in \Sigma_j$ and that $\delta_i$ converges to zero. Hence, for sufficiently large $i$, $T_i \in \mathfrak{a}(B(p,\alpha/2) \cap M;1/i,\delta_i)$. 

Then Step 2 in Construction 3.10 of Pitts \cite{pitts}, with $K=\overline B(p,\alpha/4)\cap M$, gives a sequence $T_i^*\in \mathcal{Z}_n(M;\mathbb{Z}_2)$ that satisfies
\begin{itemize}
\item[(a)] $T_i^* \in \mathfrak{a}(B(p,\alpha/2) \cap M;1/i,\delta_i)$,
\item[(b)] $0 \leq {\bf M}(T_i)-{\bf M}(T_i^*)\leq 1/i,$
\item[(c)] $T_i\llcorner (M\setminus \overline B(p,\alpha/4)) =T_i^*\llcorner (M\setminus \overline B(p,\alpha/4))$,
\item[(d)] ${\bf M}(T_i^*) \leq {\bf M}(S)$ for every $S \in \mathcal{Z}_n(M;\mathbb{Z}_2)$ satisfying $\s(S-T_i) \subset \overline B(p,\alpha/4)$ and $\mathcal{F}(S-T_i^*) \leq \delta_i$,
\item[(e)] $|T_i^*|$ is stable in $B(p,\alpha/4) \cap M$,
\item[(f)] for every $q \in B(p,\alpha/4) \cap M$, there exists $r>0$ such that ${\bf M}(S) \geq {\bf M}(T_i^*)$ for every $S \in \mathcal{Z}_n(M;\mathbb{Z}_2)$ satisfying $\s(S-T_i^*) \subset B(q,r)$.
\end{itemize}

A subsequence $\{T_j^*\}$ converges to a {\it replacement} $V^*$ of $V$ (Steps 3 and 4 of Section 3.10 \cite{pitts}). Lemma 6.11 of Zhou \cite{zhou} proves $V^*=V$. By passing to a further subsequence, we can assume $T_j^*$ converges to some $T \in \mathcal{Z}_n(M;\mathbb{Z}_2)$ in the flat topology. In particular,  $\s(T) \subset \s(V)$. By the Constancy Theorem for flat chains and property (c) above, we conclude $T=\Sigma$. On the other hand, Schoen-Simon theory \cite{schoen-simon} together with  properties (b) and (e) above gives that $|T_j^*|$ converges locally graphically and smoothly to $V$ in $M\cap  B(p,\alpha/4)$. Therefore the number of sheets of $|T_j^*|$ in $M\cap  B(p,\alpha/4)$ is equal to the multiplicity $m_j$ of the component $\Sigma_j$. Since $n_j \neq m_j {\rm \, mod \,} 2$, this is in contradiction with the fact that $T_j^*$ converges to $\Sigma$ in the flat topology.
\end{proof}


\subsection{Theorem}\label{homotopy.flat.wrong} {\it Suppose $3\leq (n+1) \leq 7$. Let $$V = m_1 \cdot |\Sigma_1| + \cdots + m_N \cdot |\Sigma_N|$$ be an embedded minimal cycle, where $\{\Sigma_1, \dots, \Sigma_N\}$ is a disjoint collection, and 
$$\Sigma = n_1\cdot  \Sigma_1 + \cdots + n_N \cdot \Sigma_N,$$ $n_i \in \{0,1\}$, be a mod 2 flat cycle  with 
$n_j \neq m_j {\rm \, mod\,} 2$ for some $1\leq j\leq N$.

Let $\Phi:X \rightarrow \mathcal{Z}_n(M; {\bf M}; \mathbb{Z}_2)$, where $X$ is a subcomplex of $I(m,t)$. There exists $\varepsilon=\varepsilon(m,V,\Sigma)>0$ such that if 
for every $x \in X$, we have $${\bf F}(|\Phi(x)|, V) \leq \varepsilon {\rm \, \, and\, \,} \mathcal{F}(\Phi(x), \Sigma)\leq \varepsilon,$$
then for every $\delta>0$ there exists a homotopy $H:[0,1] \times X \rightarrow \mathcal{Z}_n(M; {\bf M}; \mathbb{Z}_2)$ between $H(0,\cdot)=\Phi$ and $H(1,\cdot)=\Phi'$ with:
\begin{itemize}
\item[(a)] ${\bf M}(H(t,x))\leq {\bf M}(\Phi(x))+\delta$ for every $(t,x)\in [0,1]\times X$,
\item [(b)] ${\bf M}(\Phi'(x))\leq {\bf M}(\Phi(x))-\varepsilon/6$ for every $x\in X$.
\end{itemize}
}
\begin{proof}

Let $r_0>0$ be sufficiently small so that, for every $p \in \s(V)$,
\begin{itemize}
\item the restriction of $v \mapsto |v-p|^2$ to $M \cap B(p,2r_0)$ is strictly convex, 
\item $\Sigma$ is transversal to $M \cap \partial B(p,s)$ for every $s< 2r_0$,
\item if $p\in \Sigma_l$, then $B(p,2r_0) \cap \Sigma_{l'}=\emptyset$ for every $1\leq l'\leq N$, $l'\neq l$.
\end{itemize}
Let $\alpha=\frac{r_0}{16 \cdot 8^{c-1}}$, where $c=(3^m)^{3^m}$, and $\varepsilon=\varepsilon(\alpha,V,\Sigma)>0$ as in Proposition \ref{flat.wrong}.

Let $0<\delta<\delta_0(M)$ be a small number so that 
$$
m\delta (1+4(3^m-1)\frac{8^c}{r_0})< \min \{\frac{\varepsilon}{3^m8}, \gamma_{iso}\},
$$
$$C(m,\frac{r_0}{8^c})\delta < \delta_0,$$
$$(C_0C(m,\frac{r_0}{8^c})+1)\delta < \min\{\delta_1, \frac{\varepsilon}{60}\},$$
where $\delta_0=\delta_0(M)$ and $C_0=C_0(M,m+1)$ are the constants of Theorem 3.10 of \cite{marques-neves-infinitely}, $\delta_1$ is as in Proposition \ref{close.implies.homotopic}, and $C(m,\frac{r_0}{8^c})$ is as in Theorem \ref{combinatorial.thm}.

Choose  $q\in \mathbb{N}$ such that  
\begin{equation}\label{mass.fine}
\sup\{{\bf M}(\Phi(x)-\Phi(y)): x,y \mbox{ in a common cell of }X(q)\}\leq \delta.
\end{equation}
Consider the map  $\phi:X(q)_0 \rightarrow  \mathcal Z_n(M;\mathbb{Z}_2)$ given by   $\phi(x)=\Phi(x)$ for every $x \in X(q)_0$,  and let $L=\sup_{x\in X(q)_0} {\bf M}(\phi(x))$.  It satisfies ${\bf f}(\phi)\leq \delta$.   

We choose $p \in \Sigma_j$ and $\{p_1,\dots,p_c\}\subset \Sigma_j$ such that $|p_i-p|=\frac{3r_0}{4\cdot 8^{i-1}}$ for every $i=1, \dots, c$.
By Proposition \ref{flat.wrong}, for every $x\in X(q)_0$ we have that $\Phi(x)$ admits an $(\varepsilon, \delta)$-deformation in $B(p_i,\alpha)\cap M$ for every $i=1, \dots,c$. Now
$$
\overline B(p_i,\alpha)\cap M \subset A(p, r_i-s_i, r_i+s_i) \cap M,
$$
with $r_i=\frac{3r_0}{4\cdot 8^{i-1}}$ and $s_i=\frac{r_0}{8^i}$. Hence the assumptions of Theorem \ref{combinatorial.thm} are satisfied with $\mathcal W=\emptyset$, $\bar \varepsilon = R= \varepsilon$, $\eta=2r_0$, and $p(x)=p$, $r_i(x)=r_i$ and $s_i(x)=s_i$ for every $x\in X(q)_0$, $i=1, \dots, c$.

In this case, we can find $k \geq q$ and a map $\phi^*:X(k)_0\rightarrow \mathcal{Z}_n(M, \mathbb{Z}_2)$ that satisfies properties (i)-(vi) of Theorem \ref{combinatorial.thm}.  Then, by Property (i) and Theorem 3.10 of \cite{marques-neves-infinitely}, there exists an Almgren extension $\Psi: [0,1] \times X \rightarrow \mathcal{Z}_n(M;{\bf M};\mathbb{Z}_2)$ of the discrete homotopy $\psi$ between $\phi$ and $\phi^*$. It satisfies
\begin{eqnarray}\label{fineness.psi}
&& \sup \{{\bf M}(\Psi(v)-\Psi(w)): v,w {\rm \, \, in \, a \, common \, cell\, of \, }([0,1]\times X)(k)\} \nonumber\\
&&\leq C_0 {\bf f}(\psi)\nonumber\\
&&\leq C_0C(m,s_c)\delta.
\end{eqnarray}
Together with property (iii) of Theorem \ref{combinatorial.thm} and inequality (\ref{mass.fine}), this gives that
$$
{\bf M}(\Psi(t,x))\leq {\bf M}(\Phi(x))+ (1+C_0C(m,s_c)+2\cdot 3^m m(1+4(3^m-1)s_c^{-1}))\delta
$$
for every $(t,x)\in [0,1]\times X$.

 The map $\Psi_0=\Psi(0,\cdot)$ is an extension of $\phi \circ {\bf n}(k+t,q+t)$.
Hence
$$
\sup\{{\bf M}(\Phi(x)-\Psi_0(x)):x\in X\}\leq (C_0C(m,s_c)+1)\delta.
$$
Since $(C_0C(m,s_c)+1)\delta < \delta_1$, then $\Psi_0$  is homotopic to $\Phi$ in the mass topology by Proposition \ref{close.implies.homotopic}. The homotopy  $\tilde H: [0,1] \times X \rightarrow \mathcal{Z}_n(M;{\bf M};\mathbb{Z}_2)$ with $\tilde H(0,\cdot)=\Phi$ and $\tilde H(1,\cdot)=\Psi_0$ can be chosen to satisfy
\begin{eqnarray*}
&&\sup \{{\bf M}(\tilde H(t,x)-\Phi(x)): t \in [0,1], x\in X\} \\
&&\hspace{2cm}\leq C_1 \sup\{{\bf M}(\Phi(x)-\Psi_0(x)):x\in X\}\\
&&\hspace{2cm}\leq C_1(C_0C(m,s_c)+1)\delta.
\end{eqnarray*}

We can concatenate $\tilde H$ to $\Psi$ to obtain a homotopy $H: [0,1] \times X \rightarrow \mathcal{Z}_n(M; {\bf M};\mathbb{Z}_2)$ with $H(0,\cdot)=\Phi$ and $H(1, \cdot)=\Psi_1=\Psi(1,\cdot)$. Note that 
$$
{\bf M}(H(t,x))\leq {\bf M}(\Phi(x))+C_2\delta,
$$
where $C_2=C_2(M,m,s_c)>0$.

Property (vi) of Theorem \ref{combinatorial.thm} implies 
$$
{\bf M}(\phi^*(x))<L-\frac{\varepsilon}{5}
$$
for every $x \in X(k)_0$, since $\mathcal W=\emptyset$. Since $\Psi_1$ is an extension of $\phi^*$,  inequality (\ref{fineness.psi}) implies that for every $x\in X$ one has
\begin{eqnarray*}
{\bf M}(\Psi_1(x))&\leq& L-\frac{\varepsilon}{5}+C_0C(m,s_c)\delta\\
&\leq&L- \frac{\varepsilon}{6},
\end{eqnarray*}
since  $C_0C(m,s_c)\delta \leq \frac{\varepsilon}{60}$.
\end{proof}

\section{The space of cycles and $k$-sweepouts}\label{space.of.cycles}

In order to do Morse theory, one needs to first understand the topology of the underlying space. In the case of the area functional, this
is the space of closed hypersurfaces.

Note that by the Constancy Theorem for mod 2 flat chains, the boundary map
$$
\partial : {\bf I}_{n+1}(M;\Z_2) \rightarrow \mathcal Z_n(M;\Z_2)
$$
is a 2-cover. Namely, if $\partial U=\partial V$ then either $U=V$ or $U=M-V$. In the next theorem, we will prove that ${\bf I}_{n+1}(M;\Z_2)$ is contractible and that $\partial$ satisfies the lifting property.

\subsection{Theorem}\label{weak.homotopy}
{\em The space of cycles $\mathcal Z_n(M;\Z_2)$  is weakly homotopically equivalent to $\mathbb{RP}^\infty$.}

\begin{proof}
Let $f:M \rightarrow \mathbb{R}$ be a Morse function, with $f(M)=[0,1]$. We claim that the map $\Phi:\mathbb{RP}^\infty \rightarrow \mathcal Z_n(M;\Z_2)$ given by
$$
\hat\Phi([a_0:a_1:\cdots:a_k:0:\cdots]) = \partial \{x\in M: a_0+a_1f(x)+\cdots+a_kf(x)^k \leq 0\}
$$
is a weak homotopy equivalence, i.e. it induces isomorphisms in every homotopy group.  In \cite{marques-neves-infinitely} (Claim 5.5), we proved the map $\hat \Phi$ is continuous in the flat topology.

\subsection{Claim}\label{lifting} {\em Let $\Psi:I^p\rightarrow \mathcal Z_n(M;\Z_2)$ be a continuous map, $p\in \mathbb{N}$, and $U_0\in {\bf I}_{n+1}(M;\Z_2)$ be such that $\partial U_0=\Psi(0)$. Then there exists a unique continuous map $U:I^p \rightarrow  {\bf I}_{n+1}(M;\Z_2)$ such that $U(0)=U_0$ and $\partial U(x)=\Psi(x)$ for every $x\in I^p$.}
\medskip

We first prove uniqueness. Let $U, U'$ be two such maps and consider $V=U-U'$. Then $V:I^p \rightarrow  {\bf I}_{n+1}(M;\Z_2)$ satisfies $V(0)=0$ and $\partial V(x)=0$ for every $x\in I^p$. By the Constancy Theorem for mod 2 flat chains, $V(x) \in \{0,M\}$ for every $x \in I^p$. This implies the set $A=\{x\in I^p: V(x)=0\}$ is both closed and open. Since $0\in A$, we have $A=I^p$ and hence $U=U'$.

Now we prove the existence of the lifting $U$ when $p=1$. By the Isoperimetric Inequality of Federer-Fleming (see Proposition 1.11 or Corollary 1.14 of \cite{almgren}), adapted to the setting of mod 2 flat chains, there exist constants $\varepsilon_M>0$ and $\nu_M>0$ such that if $T \in \mathcal{Z}_n(M;\mathbb{Z}_2)$ satisfies $\mathcal{F}(T)<\varepsilon_M$, then there exists $W \in {\bf I}_{n+1}(M; \mathbb{Z}_2)$  with $\partial W=T$ and ${\bf M}(W) \leq \nu_M\mathcal{F}(T)$. We can choose $\varepsilon_M$ to be small and in that case $W$ is unique by the Constancy Theorem. 

 By continuity of $\Psi$, we can find a partition $0=t_0 < t_1 < \cdots < t_{q-1}<t_q=1$ such that for every $s,t \in [t_{i-1},t_i]$, $1\leq i\leq q$, we have
$$
\mathcal{F}(\Psi(s),\Psi(t))<\varepsilon_M.
$$
For $t \in [t_{i-1},t_i]$, let $W_i(t)$ be the unique element of ${\bf I}_{n+1}(M; \mathbb{Z}_2)$ with $\partial W_i(t)=\Psi(t)-\Psi(t_{i-1})$ and ${\bf M}(W_i(t))\leq \nu_M\mathcal{F}(\Psi(t)-\Psi(t_{i-1}))$. The map $W_i: [t_{i-1},t_i] \rightarrow {\bf I}_{n+1}(M; \mathbb{Z}_2)$ is continuous and $W_i(t_{i-1})=0$. For $t \in [0,t_1]$, we define $U(t)=U_0+W_1(t)$. Then $U:[0,t_1] \rightarrow {\bf I}_{n+1}(M; \mathbb{Z}_2)$ is continuous, with $U(0)=U_0$ and $\partial U(t)=\Psi(t)$ for all $t \in [0,t_1]$. Suppose that we have found a continuous map $U:[0,t_{i-1}] \rightarrow {\bf I}_{n+1}(M; \mathbb{Z}_2)$ with $U(0)=U_0$ and $\partial U(t)=\Psi(t)$ for all $t \in [0,t_{i-1}]$. Then we extend it to $[0,t_i]$ by putting $U(t)=U(t_{i-1})+W_i(t)$ for $t \in [t_{i-1}, t_i]$. The existence of the lifting $U:[0,1]  \rightarrow {\bf I}_{n+1}(M; \mathbb{Z}_2)$ follows by induction.

Suppose now $p>1$. Given $x\in I^p$, we choose a continuous path $\sigma:[0,1]\rightarrow I^p$ with $\sigma(0)=0$ and $\sigma(1)=x$ and define $\Psi_\sigma=\Psi\circ \sigma :[0,1] \rightarrow \mathcal{Z}_n(M;\mathbb{Z}_2)$. We know there exists a continuous map $U_\sigma:[0,1]  \rightarrow {\bf I}_{n+1}(M; \mathbb{Z}_2)$ with $U_\sigma(0)=U_0$ and $\partial U_\sigma(t)=\Psi_\sigma(t)$ for all $t \in [0,1]$. Then we put $U(x)=U_\sigma(1)$. Note that $\partial U(x)=\Psi(x)$.  Because $I^p$ is simply-connected, a standard argument gives that $U(x)$ does not depend on $\sigma$ and the obtained extension $U:I^p  \rightarrow {\bf I}_{n+1}(M; \mathbb{Z}_2)$ is a continuous map. This finishes the proof of the claim.

\subsection{Claim}\label{contractibility} {\em The space ${\bf I}_{n+1}(M;\mathbb{Z}_2)$ is contractible.}
\medskip

We define $H:[0,1] \times {\bf I}_{n+1}(M;\mathbb{Z}_2) \rightarrow {\bf I}_{n+1}(M;\mathbb{Z}_2)$ by
$$
H(t,U) = U \llcorner \{f\leq t\}.
$$
The map $H$ is continuous, $H(1,U)=U$ and $H(0,U)=0$ for every $U \in {\bf I}_{n+1}(M;\mathbb{Z}_2)$. This proves the claim.
\medskip

Recall that
\begin{eqnarray*}
 &&\pi_1(\mathbb{RP}^\infty,1)=\mathbb{Z}_2, {\rm \, \, and}\\
  &&\pi_k(\mathbb{RP}^\infty,1)=0 {\rm \, \, for \, \, every\, \, }  k\geq 2.
\end{eqnarray*}

If $k\geq 2$, and $\Psi:I^k \rightarrow \mathcal{Z}_n(M;\mathbb{Z}_2)$ is a continuous map with $\Psi(\partial I^k)=\{0\}$, Claim \ref{lifting} implies there exists a unique continuous map $U:I^k \rightarrow {\bf I}_{n+1}(M;\mathbb{Z}_2)$ with $U(0)=0$ and $\partial U(x)=\Psi(x)$ for every $x \in I^k$. By the uniqueness of the liftings of maps defined on $I^{k-1}$, we have $U(\partial I^k)=0$. Claim \ref{contractibility} implies $U$ is homotopically constant relative to $\partial I^k$. Hence $\Psi$ is homotopically constant relative to $\partial I^k$. This proves 
$$
\pi_k(\mathcal{Z}_n(M;\mathbb{Z}_2),0)=0
$$
for every $k\geq 2$.

Now let $\sigma:[0,1] \rightarrow \mathcal{Z}_n(M;\mathbb{Z}_2)$ be a continuous map with $\sigma(0)=\sigma(1)=0$. Claim \ref{lifting} gives a unique continuous map $U:[0,1] \rightarrow {\bf I}_{n+1}(M;\mathbb{Z}_2)$ with $U(0)=0$ and $\partial U(t)=\sigma(t)$ for every $t \in [0,1]$. Then $\partial U(1)=\sigma(1)=0$, hence either $U(1)=0$ or $U(1)=M$. If $U(1)=0$, the map $\sigma$ is homotopically constant relative to $\{0,1\}$ as in the higher-dimensional argument. Conversely, if such a homotopy exists then we can lift it and a standard argument implies $U(1)=0$. These arguments give that
$$
\pi_1(\mathcal{Z}_n(M;\mathbb{Z}_2),0)=\mathbb{Z}_2,
$$
and $\sigma: [0,1] \rightarrow \mathcal{Z}_n(M;\mathbb{Z}_2)$ given by $\sigma(t)= \partial \{f\leq t\}$ is a generator.
Since $\hat\Phi([\cos(\pi t): \sin(\pi t):0: \cdots])=\partial \{f\leq -\cot(\pi t)\},$ the map
$$
\hat\Phi_*: \pi_1(\mathbb{RP}^\infty,1) \rightarrow \pi_1(\mathcal{Z}_n(M;\mathbb{Z}_2),0)
$$
is an isomorphism. The higher homotopy groups of both spaces are trivial, thus $\hat\Phi$ is a weak homotopy equivalence.

\end{proof}

\subsection{Remark} Another way of calculating the homotopy groups of the space of cycles is through Almgren's Isomorphism Formula (\cite{almgren}):
$$
\pi_k(\mathcal{Z}_l(M,\mathbb{Z}_2),\{0\})= H_{k+l}(M, \mathbb{Z}_2),
$$
for any $k$ and $l$, but the argument is more complicated. In \cite{almgren}, Almgren used integer coefficients but the same proof applies for mod 2 coefficients.
\medskip

Theorem \ref{weak.homotopy} implies $$H^1(\mathcal Z_n(M;\Z_2);\Z_2)=\Z_2=\{0,\bar\lambda\}.$$ We call $\overline{\lambda}$ the {\it fundamental cohomology class}. It has geometric meaning, namely, if $\sigma:S^1 \rightarrow \mathcal{Z}_n(M;\mathbb{Z}_2)$ is a free loop then
$\overline \lambda \cdot [\sigma]=1$ if and only if $\sigma$ is homotopically nontrivial (here $[\sigma]$ denotes the homology class induced by $\sigma$).

\subsection{Definition} 
Let $k \in \mathbb{N}$. A continuous map $\Phi:X \rightarrow \mathcal Z_n(M;{\bf F};\Z_2)$ is called a {\it $k$-sweepout} if $\lambda = \Phi^*(\bar \lambda) \in H^1(X, \mathbb{Z}_2)$ satisfies
$$
\lambda^k = \lambda \smile \cdots \smile \lambda \neq 0 \in H^k(X,\mathbb{Z}_2),
$$
where $\smile$ denotes the cup product. We denote by $\mathcal{P}_k$ the set of all $k$-sweepouts $\Phi$ (the parameter space $X={\rm dmn}(\Phi)$ is allowed to depend on $\Phi$).

Notice that if $\Phi:X \rightarrow \mathcal Z_n(M;{\bf F};\Z_2)$ is a $k$-sweepout, then every ${\bf F}$-continuous map $\Phi'$ that is homotopic to $\Phi$ in the flat topology  is also a $k$-sweepout.

\subsection{Definition}
The {\it $k$-width} of $(M,g)$ is the number
$$\omega_k(M,g):=\inf_{\Phi\in\mathcal{P}_k}\sup_{x\in {\rm dmn}(\Phi)}{\bf M}(\Phi(x)).$$

\subsection{Remark} Because of Proposition \ref{flat.approximation}, the above definition of $k$-width coincides with the definition of $k$-width of \cite{marques-neves-infinitely} (Section 4.3). In \cite{marques-neves-infinitely}, as well as in \cite{liokumovich-marques-neves}, \cite{irie-marques-neves}, \cite{marques-neves-song}, continuity in the ${\bf F}$-metric in the definition of a $k$-sweepout is replaced by continuity in the flat topology together with the no concentration of mass property.

\section{White's local min-max theorem}\label{white.section}

In this section we  prove a slightly more general version of the local Min-max Theorem  in White \cite[Theorem 4]{white-minmax}. Similar results have been obtained in  \cite{morgan-ros} and \cite{inauen-marchese}.

\subsection{Local Min-max Theorem}\label{local.minmax}{\em Let $\Sigma$ be a closed, smooth, embedded minimal hypersurface with Morse index $k$ and multiplicity one. For every $\beta>0$, there is $\varepsilon_0>0$ and a smooth family $
\{F_v\}_{v\in \bar B^k}\subset \text{Diff}(M)$ such that 
\begin{itemize}
\item[(i)] $F_0=\text{Id}$, $F_{-v}=F_v^{-1}$ for all $v\in \bar B^k$;
\item[(ii)]  the function  $$A^{\Sigma}:\bar B^k\rightarrow [0,\infty], \quad A^{\Sigma}(v)=||(F_v)_{\#}\Sigma||(M),$$
 is strictly concave;
 \item[(iii)] $||F_v-\text{Id}||_{C^1}<\beta$ for all $v\in \bar B^k$;
\end{itemize}
and such that for every $S\in \mathcal Z_{n}(M;\Z_2)$ with $\mathcal F(S,\Sigma)<\varepsilon_0$, we have
$$\max_{v\in \bar B^k}||(F_v)_{\#}S||(M)\geq {\bf M}(\Sigma)$$
with equality only if $\Sigma=(F_v)_{\#}S$ for some $v\in \bar B^k$.
}

\begin{proof}
Let $e_1, \dots, e_k$ denote the coordinate vectors in $\R^k$. From the fact that $\Sigma$ has Morse index $k$ we can choose a smooth family $
\{F_v\}_{v\in \bar B_k}\subset \text{Diff}(M)$ with   properties (i) and (ii)  such  that 
$$\{X_i=\frac{d}{dt} {F_{te_i}}_{|t=0}\}_{i=1,\dots,k}$$
 is, on $\Sigma$,  an $L^2$-orthonormal  set of normal eigensections of the Jacobi operator. 
 
 Choose $\{\eta_i\}_{i=1}^k\subset C^{\infty}(M)$  so that $|\eta_i|\leq 1$ and such that, along $\Sigma$, $\eta_i=0$ and $\nabla \eta_i=X_i$. Given $S\in \mathcal Z_{n}(M;\Z_2)$, consider the smooth map
$$P^S:\bar B^k\rightarrow \R^k, \quad P^S(v)=\sum_{i=1}^k||(F_v)_{\#}S||(\eta_i)\cdot e_i.$$
Note that $DP^{\Sigma}(0)=\text{Id}$. Therefore there exists $\delta>0$ such that  $P^{\Sigma}_{|\overline B^k_\delta(0)}$ is a diffeomorphism onto its image, with $P^\Sigma(0)=0$. We can choose $\delta>0$ sufficiently small so that $||F_v-\text{Id}||_{C^1}<\beta$ for all $v\in \bar B^k_\delta(0).$

Let $\lambda_1$ be the first eigenvalue of the Jacobi operator of $\Sigma$, and  $c_0=|\lambda_1|+1.$  Consider the functional on $\mathcal Z_{n}(M;\Z_2)$ given by
$${A}^*(S)=||S||(M)+c_0\sum_{i=1}^k \left(||S||(\eta_i)\right)^2.$$
One can check that $\Sigma$ is  a strictly stable critical point for $A^*$ (see \cite[page 215]{white-minmax}).

The functional $A^*$ is continuous in the ${\bf F}$-topology and lower semicontinuous under flat convergence, as proven by White  in \cite[Claim 2, page 212]{white-minmax}.  

\subsection{Proposition}\label{prop.local.minmax}{\em There exists $\varepsilon_1>0$ so that for every $S\in \mathcal Z_{n}(M;\Z_2)$ with $\mathcal F(S,\Sigma)<\varepsilon_1$,
we have ${A}^*(S)> {A}^*(\Sigma)$ unless $S=\Sigma$.
}

\begin{proof}
We argue by contradiction and assume there is a sequence $\{S_j\}$ converging to $\Sigma$ in the flat topology such that ${A}^*(S_j)\leq  {A}^*(\Sigma)={\bf M}(\Sigma)$  and $S_j\neq \Sigma$ for all $j\in \N$.  One immediate consequence is that $||S_j||(M)\leq {\bf M}(\Sigma)$ for all $j\in \N$. Lower semicontinuity of the mass functional under flat convergence implies $\lim_{j\rightarrow \infty} ||S_j||(M)={\bf M}(\Sigma)$.  Hence  $\{S_j\}$  also converges to $\Sigma$ in the ${\bf F}$-topology  (\cite[2.1 (20)]{pitts}).

From Theorem 5 in \cite{white-minmax} we obtain the existence of a tubular neighborhood $U$ of $\Sigma$ so that $\Sigma$ is a strict minimizer of $A^*$ among all mod $2$ cycles contained in $U$ and homologous to $\Sigma$. Thus $S_j \cap (M\setminus U)\neq \emptyset $ for all  $j\in \N$. Choose a smaller tubular neighborhood $V$ with closure contained in $U$ and set
$$m_j=\inf\{A^*(T): \text{spt}\,(T-S_j)\subset M\setminus V, T\in  \mathcal Z_{n}(M;\Z_2)\}.$$
From lower semicontinuity and the Federer-Fleming Compactness Theorem we obtain $T_j\in  \mathcal Z_{n}(M;\Z_2)$ so that $A^*(T_j)=m_j$ and $ \text{spt}\,(T_j-S_j)\subset M \setminus V$.
\medskip

\noindent{\bf Claim 1:} {\em $\{T_j\}$ converges to $\Sigma$ in the ${\bf F}$-topology.}
\medskip

We have that $||S_j||(M\setminus V)\to 0$ as $j \to \infty$. We have $||T_j||(V)=||S_j||(V)$ and $||S_j||(V)- ||S_j||(M)\to 0$. On the other hand,
$$||T_j||(M)\leq A^*(T_j)\leq A^*(S_j)=||S_j||(M) +c_0\sum_{i=1}^k \left(||S_j||(\eta_i)\right)^2.$$
and hence $||T_j||(M\setminus V)\to 0$. Thus ${\bf M}(S_j-T_j)\to 0$, proving the claim.  

\medskip

\noindent{\bf Claim 2:} {\em For every  sufficiently large $j$,  $\s(T_j)\subset U$.}

\medskip

Denote the first variation of  $|T_j|$  with respect to the $n$-dimensional area function  by $\delta |T_j|$. Since $T_j$ is a critical point for $A^*$ in $M \setminus \bar V$, for any vector field $X$ compactly supported in $M\setminus \bar V$ we have   $$\delta |T_j|(X)+2c_0\sum_{i=1}^k||T_ j||(\eta_i)\delta |T_j| (\eta_i X)+2c_0\sum_{i=1}^k||T_ j||(\eta_i)||T_j|| (\langle X^{\bot}, \nabla \eta_i\rangle)=0,$$
where $X^{\bot}$ denotes the  projection of $X$ onto the orthogonal complement of $T_x T_j$, which is  defined $\mathcal{H}^n$ a.e. on $T_j$. Hence, for every sufficiently large $j$ and any open set $B\subset M\setminus \bar V$, 
$$||(\delta |T_j|)||(B)\leq a_j||T_j||(B)$$
for some sequence $a_j\to 0$. 

This gives that $T_j$ satisfies a monotonicity formula in $M\setminus \bar V$ with uniform constants (see Section 40 of \cite{simon}).
It implies that if $$\s(T_j) \cap (M \setminus U)\neq \emptyset,$$ then there is a uniform $a>0$ such that $||T_j||(M \setminus V) \geq a$. Contradiction, hence Claim 2 is proved.
\medskip

Since $A^*(T_j) \leq A^*(\Sigma)$,  we deduce from Theorem 5 \cite{white-minmax} that $T_j=\Sigma$ for all $j$ large enough. Because  $\text{spt}\,(T_j-S_j)\subset M\setminus V$, we get $S_j=\Sigma + S_j \llcorner (M \setminus V)$. Now since $||S_j||(M)\leq {\bf M}(\Sigma)$, we deduce that $S_j=\Sigma$ for all  $j$ large enough. This gives a contradiction, which proves the proposition.

\end{proof}

We can choose $0<\delta'<\delta$ such that $\mathcal F((F_{v})_{\#}\Sigma, \Sigma)< \varepsilon_1/2$ for every $v\in \overline B^k_{\delta'}$. To finish the proof of the Local Min-max Theorem, we argue again by contradiction and assume there is a sequence $\{S_j\}$ converging to $\Sigma$ in the flat topology such that 
$(F_v)_{\#}S_j\neq \Sigma$ for all $v\in \bar B^k$ and
$$\max_{v\in \bar B^k_{\delta'}} ||(F_v)_{\#}S_j||(M)\leq ||\Sigma||(M).$$
This implies $||S_j||(M) \leq ||\Sigma||(M)$, so  as before we deduce that $S_j$ converges to $\Sigma$ in the ${\bf F}$-topology. This implies uniform convergence of  $P^{S_j}_{|\partial \bar B^k_{\delta'}}$  to  $P^\Sigma_{|\partial \bar B^k_{\delta'}}$. A degree argument gives the existence of $v_j\in  B^k_{\delta'}$ such that   $P^{S_j}(v_j)=0$.
Hence, for every  sufficiently large $j$, we obtain from  Proposition \ref{prop.local.minmax} that
$$||(F_{v_j})_{\#}S_j||(M)=A^*((F_{v_j})_{\#}S_j)>A^*(\Sigma)=||\Sigma||(M).$$
Contradiction, and the theorem is proved.

\end{proof}

\section{Morse Index of Multiplicity One Min-max Minimal Hypersurfaces}\label{proof.of.main.theorem}

The next proposition uses the index estimates of \cite{marques-neves-index}.

\subsection{Proposition}\label{find.homotopy.class}
{\it Suppose $(M^{n+1},g)$ is a bumpy metric, $3\leq (n+1) \leq 7$. Then, for each $k \in \mathbb{N}$, there exists a homotopy class $\Pi$ of $k$-sweepouts such that
$$
\omega_k(M,g)={\bf L}(\Pi).
$$
}

\begin{proof}
Fix $k\in \mathbb{N}$. Choose a sequence $\{\Phi_i\}\subset  \mathcal{P}_k$ such that
$$
\lim_{i \rightarrow \infty} \sup\{M(\Phi_i(x)): x\in X_i={\rm dmn}(\Phi_i)\} = \omega_k(M,g).
$$

Denote by $X_i^{(k)}$ the $k$-dimensional skeleton of $X_i$. Then $H^k(X_i,X_i^{(k)};\Z_2)=0$ and hence the long exact cohomology sequence gives that the natural pullback map from $H^k(X_i;\Z_2)$ into $H^k(X_i^{(k)};\Z_2)$ is injective. This implies  ${(\Phi_i)}_{|X_i^{(k)}}\in \mathcal{P}_k$. The definition of $\omega_k$ then implies 
$$
\lim_{i \rightarrow \infty} \sup\{M(\Phi_i(x)): x\in X_i^{(k)}\} = \omega_k(M,g).
$$

Let $\Pi_i$ denote the homotopy class of $(\Phi_i)_{|X_i^{(k)}}$.
Then
$$
 \omega_k(M,g) \leq {\bf L}(\Pi_i) \leq  \sup\{M(\Phi_i(x)): x\in X_i^{(k)}\},
$$
and in particular $$\lim_{i \rightarrow \infty} {\bf L}(\Pi_i) = \omega_k(M,g).$$

Theorem 1.2 of \cite{marques-neves-index}   implies the existence of a finite disjoint collection $\{\Sigma_{i,1},\dots,\Sigma_{i,N_i}\}$ of closed, smooth, embedded minimal hypersurfaces in $M$, and integers $\{m_{i,1},\dots,m_{i,N_i}\} \subset \mathbb{N}$, such that
$$
 {\bf L}(\Pi_i) = \sum_{j=1}^{N_i} m_{i,j} ||\Sigma_{i,j}||(M),
$$
and
$$
\sum_{j=1}^{N_i} {\rm index}(\Sigma_{i,j}) \leq k.
$$

Since the metric $g$ is bumpy, the Compactness Theorem of Sharp (Theorem 2.3 of \cite{sharp},  Remark 2.4) implies there are finitely many minimal  hypersurfaces with area bounded by $\omega_k(M,g)+1$ and index bounded by $k$. This gives that $\{{\bf L}(\Pi_i)\}$ is a finite set, so there exists $i_0\in \mathbb{N}$ with $\omega_k(M,g) = {\bf L}(\Pi_{i_0})$.

\end{proof}

The main theorem of this paper is:

\subsection{Theorem}\label{deformation.theorem}
{\it Suppose $(M^{n+1},g)$ is a bumpy metric, $3\leq (n+1) \leq 7$, and let $\Pi$ be a homotopy class of $k$-sweepouts with ${\bf L}(\Pi)=\omega_k(M,g)$. Suppose $\{\Phi_i\}$ is a minimizing sequence in $\Pi$ such that every embedded minimal cycle of ${\bf C}(\{\Phi_i\})$ has multiplicity one. Then
there exists an embedded minimal cycle  $\Sigma \in {\bf C}(\{\Phi_i\})$ with
$$
{\rm index}(\Sigma) = k.
$$
}

\begin{proof}
Let $X={\rm dmn}(\Phi_i)$. We will think of $X$ as a simplicial complex in this proof.

Let $L={\bf L}(\Pi)$ and $\mathcal W_L$ be the set of all  stationary integral varifolds  in $M$ with mass equal to $L$  and whose support is a closed, smooth, embedded, minimal hypersurface (equivalently, the set of all embedded minimal cycles in $M$ with mass equal to $L$).
We denote by $\mathcal W_{L,j}$, $\mathcal W_{L}^j$ the elements in $\mathcal W_L$ whose support has Morse index less than or equal to $j$  and bigger than or equal to $j$, respectively. Sharp's  Compactness Theorem \cite{sharp} implies $\mathcal W_{L,j}$ is a finite set while $\mathcal W_{L}^j$ is countable, because the metric $(M,g)$ is bumpy.

Let $\alpha>0$ be such that every embedded minimal cycle in  $\overline {\bf B}^{\bf F}_\alpha({\bf C}(\{\Phi_i\})) \cap \mathcal{W}_{L,k}$ has multiplicity one. The existence of $\alpha$ is again a consequence of Sharp's Compactness Theorem. By successively applying Deformation Theorem A of  Section 5 of \cite{marques-neves-index}  for each $\tilde{\Sigma} \in \mathcal{W}_L^{k+1}$ and $$K = \{V \in \mathcal{V}_n(M): ||V||(M) \leq L+1, {\bf F}(V, {\bf C}(\{\Phi_i\}))\geq \alpha\},$$ we  construct a minimizing sequence $\{\Psi_i\}\subset \Pi$  such that ${\bf C}(\{\Psi_i\}_{i\in\N})\cap \mathcal W_{L}^{k+1}=\emptyset$ and 
$$
{\bf C}(\{\Psi_i\}) \subset \overline {\bf B}^{\bf F}_\alpha({\bf C}(\{\Phi_i\})).
$$
We can assume that it is pulled tight.  Note that every embedded minimal cycle in  ${\bf C}(\{\Psi_i\})$ has multiplicity one.  

Set $\Lambda_1={\bf C}(\{\Psi_i\})\cap \mathcal W_{L,k}=\{\tilde{\Sigma}_1, \dots, \tilde{\Sigma}_q\}$. Let $[\tilde{\Sigma}_i]$ be the mod 2  flat cycle (not necessarily a boundary) induced by $\tilde{\Sigma}_i$, so that the density function satisfies $\Theta([\tilde{\Sigma}_i],x)= \Theta(\tilde{\Sigma}_i,x) {\rm \, mod\,}2$ for every $x\in M$. We would like to show that there exists an element $\tilde{\Sigma} \in \Lambda_1$ with ${\rm index}(\tilde{\Sigma})=k$. Suppose not, by contradiction.

Let 
$$
\varepsilon_i=\min\{\varepsilon(m,\tilde{\Sigma}_i,T): T \in \mathcal Z_n(M;\Z_2), \s(T)\subset \tilde{\Sigma}_i, T \neq [\tilde{\Sigma}_i]\},
$$
and $\varepsilon =\min\{\varepsilon_1, \dots, \varepsilon_q\}$, where $\varepsilon(m,V,\Sigma)$ is as in Theorem \ref{homotopy.flat.wrong}.
Consider the finite set $\mathcal T$ of all elements in $\mathcal Z_n(M;\Z_2)$ whose support is contained in the support of some element in $\Lambda_1$. We choose 
$$0<s_1<\min\{\varepsilon/(12),\varepsilon_0(\tilde{\Sigma}_{1}), \dots, \varepsilon_0(\tilde{\Sigma}_{q}), \delta(M,m)/8\},$$ where $\varepsilon_0$ is the constant of Theorem \ref{local.minmax} and $\delta(M,m)$ is the constant of Proposition 3.5 of \cite{marques-neves-infinitely}, so that  for all $\Sigma\in \Lambda_1$, $\Sigma'\in \mathcal T$,
$${\bf B}_{2s_1}^{\bf F}(\Sigma)\cap \Lambda_1 =\{\Sigma\}, \quad {\bf B}_{2s_1}^{\mathcal F}(\Sigma')\cap \mathcal T =\{\Sigma'\}.$$
There exists $0<s_2<s_1$ so  that for all $T\in \mathcal Z_n(M;\Z_2)$ with ${\bf F}(|T|,\Sigma)<2s_2$ for some $\Sigma\in \Lambda_1$ there is $\Sigma'\in \mathcal T$ so that ${\mathcal F}(T,\Sigma')<s_1$. This is a consequence of the Constancy Theorem of flat chains and lower semicontinuity of mass with respect to flat topology. 

Given $0<s<s_2$, we apply Theorem \ref{theorem.improv.minmax} to obtain  $\{\hat \Phi_i\}\subset \Pi$ with the following properties:
\begin{itemize}
\item[(i)] ${\bf C}(\{ \hat \Phi_i\})\cap \mathcal W^{k+1}_L= \emptyset$;
\item[(ii)] for some $0<\gamma<L$ and for all $i$ sufficiently large
$$\{|\hat\Phi_i(x)|\in\mathcal V_n(M):{\bf M}(\hat \Phi_i(x))\geq L-\gamma\}\subset \bigcup_{\Sigma\in \Lambda_1}{\bf B}_s^{\bf F}(\Sigma).$$
\end{itemize}
We can suppose $\hat \Phi_i$ is continuous in the mass norm, by Proposition \ref{M.regularization}.

Fix $i$ sufficiently large, and let $0<\delta<\min\{s, \gamma/4\}$. Denote by $X^k$ the $k$-dimensional skeleton of $X$. Then ${(\hat\Phi_i)}_{|X^k}$ is also a $k$-sweepout (see proof of Proposition \ref{find.homotopy.class}). Now let $\tilde{X}^k$ be the union of all $k$-dimensional simplices of $X^k$. If $Z$ is the union of all simplices in $X^k$ that are not faces of a $k$-dimensional simplex, then an application of the Mayer-Vietoris sequence of $X^k = \tilde{X}^k \cup Z$ gives that the restriction map
$$
r^*: H^k(X^k,\mathbb{Z}_2) \rightarrow H^k(\tilde{X}^k,\mathbb{Z}_2)
$$
is an isomorphism. Hence ${(\hat\Phi_i)}_{|\tilde{X}^k}$ is also a $k$-sweepout.

Since cohomology with $\mathbb{Z}_2$ coefficients is dual to homology with $\mathbb{Z}_2$ coefficients, the property that
$\hat\Phi_i^*(\overline{\lambda}^k)_{|\tilde{X}^k} \neq 0 \in H^k(\tilde{X}^k, \mathbb{Z}_2)$ means there exists $\sigma \in H_k(\tilde{X}^k,\mathbb{Z}_2)$ such that
$$
\hat\Phi_i^*(\overline{\lambda}^k)_{|\tilde{X}^k} \cdot \sigma =1.
$$
The equivalence between simplicial homology and singular homology implies there exist $k$-dimensional simplices $t_1, \dots, t_h$ in 
$\tilde{X}^k$ with $\sigma = [\sum_{l=1}^h t_l]$. In particular, $\sum_{l=1}^h \partial t_l=0$. Therefore
\begin{equation}\label{homology.condition}
\overline{\lambda}^k \cdot \left[\sum_{l=1}^h (\hat\Phi_i)_{\#}(t_l)\right] = 1,
\end{equation}
where $(\hat\Phi_i)_{\#}(t_l)$ denotes the singular simplex associated with $(\hat{\Phi}_i)_{|t_l}$. Notice that the choice of $\{t_1, \dots, t_h\}$ depends on $i$. We will omit this dependence in what follows.

Let $Y_i$ be the $m_i$'th successive barycentric subdivision of $\cup_{l=1}^h t_l$ so that ${\bf F}(\hat\Phi_i(x),\hat\Phi_i(y))<\delta$ whenever $x,y$ belong to a common simplex in  $\cup_{l=1}^h t_l$. Let $W_i$ be the union of all $k$-dimensional simplices $t \in Y_i$ such that ${\bf M}(\hat\Phi_i(x)) \geq L-\gamma/2$ for some  $x\in t$. Then ${\bf M}(\hat\Phi_i(y))\geq L-\gamma$ for every $y\in W_i$.
Let $W_{i,1}, \dots, W_{i,r}$ be the connected components of $W_i$ ($r$ depending on $i$). Then by property (ii) above and the choice of $s$ we have that, for each $1 \leq p\leq r$ there exists $1\leq q_p\leq q$ with
\begin{equation}\label{F.close}
{\bf F}(|\hat\Phi_i(y)|, \tilde{\Sigma}_{q_p}) < s
\end{equation}
for every $y \in W_{i,p}$. Then there exists a unique $T_p \in \mathcal{Z}_n(M;\mathbb{Z}_2)$ with $\s(T_p) \subset \tilde{\Sigma}_{q_p}$ and
\begin{equation}\label{flat.close}
\mathcal{F}(\hat\Phi_i(y), T_p)<s_1
\end{equation}
for every $y \in W_{i,p}$.

Notice that $\partial (\sum_{t\in Y_i}t)=0$. This means that every $(k-1)$-simplex of $Y_i$ is a face of an even number of $k$-simplices of $Y_i$. If $u \in \partial (\sum_{t\in W_i} t)$, then $u$ is a face of some $k$-simplex in $W_i$ and of some $k$-simplex in $Y_i\setminus W_i$. Hence for every $y\in u$, we have
\begin{equation}\label{area.down}
L-\gamma \leq {\bf M}(\hat\Phi_i(y))<L-\gamma/2.
\end{equation}

Let us fix one connected component $W_{i,p}$, and suppose ${\rm index}(\tilde{\Sigma}_{q_p})=j < k$. We proceed differently, depending on whether $T_p = [\tilde{\Sigma}_{q_p}]$ or $T_p \neq  [\tilde{\Sigma}_{q_p}]$. Suppose $T_p = [\tilde{\Sigma}_{q_p}]$, and let $
\{F_v\}_{v\in \bar B_j}\subset \text{Diff}(M)$ be the smooth family associated with $\tilde{\Sigma}_{q_p}$ given by Theorem \ref{local.minmax}. We can arrange so that
$$
\mathcal F((F_v)_{\#}(T_p),T_p) < \delta(M,m)/2,
$$
$$
\mathcal{F}((F_v)_{\#}(S_1), (F_v)_{\#}(S_2))\leq 2 \mathcal{F}(S_1,S_2),
$$
and
$$
{\bf F}((F_v)_{\#}(V_1),(F_v)_{\#}(V_2))\leq 2 {\bf F}(V_1,V_2)
$$
for every $v \in \overline B^j$, $S_1,S_2 \in \mathcal{Z}_n(M;\mathbb{Z}_2)$ and $V_1,V_2\in \mathcal V_n(M)$.
For $y \in {\rm support}(\partial (\sum_{t\in W_{i,p}} t))$, let $A^y:\overline{B}^j \rightarrow [0,\infty)$ be the function
$$
A^y(v) =||{(F_v)}_{\#}(|\hat\Phi_i(y)|)||(M).
$$
If $s$ is sufficiently small, the function $A^y$ is strictly concave and has a unique maximum at $m(y) \in B^m_{\frac12}(0)$. The function
$$m: {\rm support}(\partial (\sum_{t\in W_{i,p}} t)) \rightarrow B^m_{\frac12}(0)$$
is continuous. Because of the upper bound in (\ref{area.down}), Theorem \ref{local.minmax} applied with $S=\hat\Phi_i(y)$ gives that $m(y) \neq 0$ for every $y \in {\rm support}(\partial (\sum_{t\in W_{i,p}} t))$.  Hence we can find $\eta>0$ such that $\eta \leq |m(y)| <1/2$ for every $y \in {\rm support}(\partial (\sum_{t\in W_{i,p}} t))$.

Consider the one-parameter flow $\{\phi^y(\cdot, t)\}_{t \geq 0} \subset {\rm Diff}(\overline{B}^j)$ generated by the vector field
$$
v \mapsto -(1-|v|^2)\nabla A^y(v).
$$
With $v\in \overline{B}^j$ fixed, the function $t \mapsto A^y(\phi^y(v,t))$ is non-decreasing. It is strictly decreasing except if $v=m(y)$ or if 
$v \in \partial B^j$, in which cases it is constant. We also have that
$$
\lim_{t\rightarrow \infty}\phi^y(v,t) \in \partial B^j
$$
if $v \neq m(y)$, and the limit is uniform if $|v-m(y)| \geq \eta$.

Therefore there exists $t_0>0$ such that the homotopy
$$
H_1: [0,1] \times {\rm support}(\partial (\sum_{t\in W_{i,p}} t)) \rightarrow \mathcal{Z}_n(M;{\bf F}; \mathbb{Z}_2)
$$
given by
$$
H_1(t,y) = {(F_{\phi^y(0,t_0t)})}_{\#}(\hat\Phi_i(y))
$$
satisfies
\begin{itemize}
\item $H_1(0,y) = \hat\Phi_i(y)$,
\item ${\bf F}(H_1(1,y), {(F_{w(y)})}_{\#}(\hat\Phi_i(y)))< \delta$ for some continuous function $$w: {\rm support}(\partial (\sum_{t\in W_{i,p}} t))  \rightarrow \partial B^j,$$
\item ${\bf M}(H_1(t,y)) < L-\gamma/2$ 
\end{itemize}
for every $(t,y) \in [0,1] \times {\rm support}(\partial (\sum_{t\in W_{i,p}} t))$. Note that
\begin{eqnarray*}
&&\mathcal F(H_1(t,y),T_p) \leq \mathcal F({(F_{\phi^y(0,t_0t)})}_{\#}(\hat\Phi_i(y)), {(F_{\phi^y(0,t_0t)})}_{\#}(T_p))\\
&&\hspace{3cm}+ \mathcal F ({(F_{\phi^y(0,t_0t)})}_{\#}(T_p), T_p)\\
&&\leq 2 \mathcal F(\hat\Phi_i(y), T_p)+\delta(M,m)/2\\
&&\leq 2s_1 + \delta(M,m)/2\\
&&\leq 3\delta(M,m)/4
\end{eqnarray*}
for every $(t,y) \in [0,1] \times {\rm support}(\partial (\sum_{t\in W_{i,p}} t))$.

Because of inequality (\ref{F.close}), we have
\begin{eqnarray*}
&& {\bf F}(H_1(1,y), (F_{w(y)})_{\#}(\tilde{\Sigma}_{q_p}))< \delta+{\bf F}({(F_{w(y)})}_{\#}(\hat\Phi_i(y)), (F_{w(y)})_{\#}(\tilde{\Sigma}_{q_p}))\\
 && \leq \delta+ 2s_1+2s\leq 5s_1,
 \end{eqnarray*}
 for every $y \in {\rm support}(\partial (\sum_{t\in W_{i,p}} t))$.
Note that ${\bf M}((F_{w})_{\#}(\tilde{\Sigma}_{q_p})) \leq L-b$ for every $w\in \partial B^j$ and some $0<b<\min\{\delta(M,m),\gamma\}$. 

By applying Theorem \ref{homotopy.F.thm} with the compact set $$\mathcal K=\{(F_w)_{\#}(\tilde{\Sigma}_{q_p}): w\in \partial B^j\}\subset \mathcal{Z}_n(M;{\bf F};\mathbb{Z}_2),$$ 
we conclude that if $s_1$ is sufficiently small (depending on $\mathcal{K}$ and $b$),  there 
exists a second homotopy
$$
H_2: [0,1] \times {\rm support}(\partial (\sum_{t\in W_{i,p}} t)) \rightarrow \mathcal{Z}_n(M;{\bf F}; \mathbb{Z}_2)
$$
such that
\begin{itemize}
\item $H_2(0,y) = H_1(1,y)$,
\item $H_2(1,y)=(F_{w(y)})_{\#}(\tilde{\Sigma}_{q_p})$,
\item ${\bf F}(H_2(t,y),(F_{w(y)})_{\#}(\tilde{\Sigma}_{q_p})) \leq b/4$
\end{itemize}
for every $(t,y) \in [0,1] \times {\rm support}(\partial (\sum_{t\in W_{i,p}} t))$. In particular,
$$
{\bf M}(H_2(t,y)) \leq {\bf M}((F_{w(y)})_{\#}(\tilde{\Sigma}_{q_p}))+b/4< L-b/2,
$$
and
\begin{eqnarray*}
&&\mathcal F(H_2(t,y),T_p) \leq \mathcal F(H_2(t,y),(F_{w(y)})_{\#}(T_p)) + \mathcal F((F_{w(y)})_{\#}(T_p),T_p)\\
&&\leq b/4+\delta(M,m)/2\\
&&\leq 3\delta(M,m)/4
\end{eqnarray*}
for every $(t,y) \in [0,1] \times {\rm support}(\partial (\sum_{t\in W_{i,p}} t))$.

If $j=0$,  Theorem \ref{local.minmax} implies ${\bf M}(\hat\Phi_i(y)) \geq L={\bf M}(\tilde{\Sigma}_{q_p})$ for every $y \in W_{i,p}$. Hence (\ref{area.down}) can be true only if $\partial (\sum_{t\in W_{i,p}} t)=0$.

Let $W_i^{(0)}$ be the union of all components $W_{i,p}$ such that ${\rm index}(\tilde{\Sigma}_{q_p})=0$ and $T_p = [\tilde{\Sigma}_{q_p}]$. Thus identity (\ref{homology.condition}) gives
\begin{eqnarray*}
1&=& \overline{\lambda}^k \cdot \left[\sum_{t \in Y_i} (\hat\Phi_i)_{\#}(t)\right]\\
&=&  \overline{\lambda}^k \cdot \left[\sum_{t \in Y_i\setminus W_i^{(0)}} (\hat\Phi_i)_{\#}(t)  + \sum_{t \in W_i^{(0)}} (\hat\Phi_i)_{\#}(t)\right]\\
&=&  \overline{\lambda}^k \cdot \left[\sum_{t \in Y_i\setminus W_i^{(0)}} (\hat\Phi_i)_{\#}(t)\right]  +  \overline{\lambda}^k \cdot \left[\sum_{t \in W_i^{(0)}} (\hat\Phi_i)_{\#}(t)\right]\\
&=&  \overline{\lambda}^k \cdot \left[\sum_{t \in Y_i\setminus W_i^{(0)}} (\hat\Phi_i)_{\#}(t)\right].
\end{eqnarray*}
The third equality above holds because $\partial(\sum_{t \in W_i^{(0)}} (\hat\Phi_i)_{\#}(t))=0$ and the fourth equality holds because
$$
\hat\Phi_i(W_i^{(0)}) \subset \bigcup_{q'=1}^q B^{\mathcal{F}}_{s_1}(\tilde{\Sigma}_{q'}).
$$

Suppose $j\geq 1$. We have  $H_{k-1}(\partial B^j, \mathbb{Z}_2)=0$, since $j<k$. This implies
$$
[w_{\#}(\partial (\sum_{t\in W_{i,p}} t))]=0 
$$
in homology. This means there exists a $k$-dimensional singular chain $\sum_j \alpha_j$ in $\partial B^j$ such that
$\sum_j \partial\alpha_j= w_{\#}(\partial (\sum_{t\in W_{i,p}} t))$. Consider the associated singular simplex $\hat{\alpha}_j:\Delta^k \rightarrow \mathcal{Z}_n(M;{\bf F}; \mathbb{Z}_2)$ given by 
$$
\hat{\alpha}_j(y) =F_{\alpha_j(y)}(\tilde{\Sigma}_{q_p}), \hspace{0.5cm} y\in \Delta^k.
$$
Here $\Delta^k$ denotes the standard $k$-simplex.

 Consider the singular chain
 \begin{eqnarray*}
&&z_{i,p}= \sum_{t\in W_{i,p}} (\hat{\Phi}_i)_{\#}(t) + (H_1)_{\#}\left([0,1] \times {\rm support}(\partial (\sum_{t\in W_{i,p}} t))\right)\\
&&+  (H_2)_{\#}\left([0,1] \times {\rm support}(\partial (\sum_{t\in W_{i,p}} t))\right) + \sum_j \hat{\alpha}_j.
 \end{eqnarray*}
 Then $\partial z_{i,p}=0$. Since
 $$
 {\rm image}(z_{i,p}) \subset B_{\delta(M,m)}^{\mathcal{F}}(\tilde{\Sigma}_{q_p}),
 $$
 we have that
 $$
 \overline{\lambda}^k \cdot [z_{i,p}]=0.
 $$
 
  In particular, the $k$-dimensional cycle
 $$
\tilde{z} = \sum_{t \in Y_i\setminus W_i^{(0)}} (\hat\Phi_i)_{\#}(t)+\sum_{p \in \mathcal{I}} z_{i,p}
 $$
 satisfies $\overline{\lambda}^k \cdot [\tilde{z}]=1$, 
 where $\mathcal{I} = \{p: {\rm index}(\tilde{\Sigma}_{q_p})>0, T_p=\tilde{\Sigma}_{q_p}\}$.

 Now suppose $T_p \neq \tilde{\Sigma}_{q_p}$. Since $s_1<\varepsilon/2$, inequalities (\ref{F.close}) and (\ref{flat.close}) give that
$$
{\bf F}(|\hat\Phi_i(y)|, \tilde{\Sigma}_{q_p}) < \varepsilon_{q_p}, \mathcal{F}(\hat\Phi_i(y), T_p)<\varepsilon_{q_p}
$$
for every $y \in W_{i,p}$. Theorem \ref{homotopy.flat.wrong}  gives a homotopy 
$$
H_3:[0,1] \times W_{i,p} \rightarrow \mathcal{Z}_n(M;{\bf M};\mathbb{Z}_2)
$$
such that 
\begin{itemize}
\item $H_3(0,y)=\hat \Phi_i(y)$,
\item ${\bf M}(H_3(t,y))\leq {\bf M}(\hat \Phi_i(y))+\delta$,
\item ${\bf M}(H_3(1,y))\leq {\bf M}(\hat \Phi_i(y))-\varepsilon/6\leq L+s-\varepsilon/6\leq L-\frac{\varepsilon}{12}$
\end{itemize}
for every $(t,y) \in [0,1] \times {\rm support}(\sum_{t\in W_{i,p}} t)$. In particular,
$$
{\bf M}(H_3(t,y))\leq L-\gamma/2+\delta\leq L-b/4
$$
 if $y \in  {\rm support}(\partial (\sum_{t\in W_{i,p}} t))$.

Consider the singular chain
 \begin{eqnarray*}
&&v_{i,p}= \sum_{t\in W_{i,p}} (\hat{\Phi}_i)_{\#}(t) +\sum_{t\in W_{i,p}} (H_3(1,\cdot))_{\#}(t)\\
&&+  (H_3)_{\#}\left([0,1] \times \partial (\sum_{t\in W_{i,p}} t)\right).
 \end{eqnarray*}
Then 
$$
v_{i,p}=\partial \left((H_3)_{\#}\left([0,1] \times  (\sum_{t\in W_{i,p}} t)\right)\right).
$$

Now the $k$-dimensional cycle
 $$
\tilde{\tilde{z}} = \sum_{t \in Y_i\setminus W_i^{(0)}} (\hat\Phi_i)_{\#}(t)+\sum_{p \in \mathcal{I}} z_{i,p}+\sum_{p\in \mathcal{J}} v_{i,p},
 $$
 where $\mathcal{J}=\{p: T_p\neq \tilde{\Sigma}_{q_p}\},$ is homologous to $\tilde{z}$. Hence $\overline{\lambda}^k \cdot [\tilde{\tilde{z}}]=1$. It follows from construction 
  that ${\bf M}(T)\leq L-\min\{b/4, \varepsilon/(12)\}$ for every $T \in {\rm image}(\tilde{\tilde{z}})$. 
  
  From $\tilde{\tilde{z}}$ we can construct a $\Delta$-complex (see Section 2.1 in \cite{hatcher}) $\tilde{\tilde{Z}}$ and a continuous map
 $\Xi_i: \tilde{\tilde{Z}} \rightarrow \mathcal{Z}_n(M;{\bf F}; \mathbb{Z}_2)$ such that
 \begin{itemize}
 \item $\Xi_i^*(\overline{\lambda})^k \neq 0 \in H^k(\tilde{\tilde{Z}}, \mathbb{Z}_2),$
 \item and ${\bf M}(\Xi_i(y)) \leq L-\min\{b/4, \varepsilon/(12)\}$ for every $y \in \tilde{\tilde{Z}}$.
 \end{itemize}
 Since every $\Delta$-complex is homeomorphic to a simplicial complex, the map $\Xi_i$ is a $k$-sweepout. Contradiction, because $L-\min\{b/4, \varepsilon/(12)\} < \omega_k(M,g)$.
 
 Hence there must exist an embedded minimal cycle $\tilde{\Sigma}$ of multiplicity one and ${\rm index}(\tilde\Sigma)=k$ in  $\overline{\bf B}^{\bf F}_\alpha({\bf C}(\{\Phi_i\})).$ Since $\alpha$ can be chosen arbitrarily small, Sharp's Compactness Theorem (\cite{sharp}) and the fact that $g$ is bumpy imply that there exists an embedded minimal cycle $\Sigma \in {\bf C}(\{\Phi_i\})$  with ${\rm index}(\Sigma)=k$ and we are done.

\end{proof}

\section{Addendum}

After this work was completed, X. Zhou \cite{zhou-multiplicity} in an exciting new work used a  novel regularization of the area functional (developed by him and Zhu in \cite{zhou-zhu-prescribed}) to prove the  Multiplicity One Conjecture  \ref{m.o.c} proposed by the authors (see also \cite{marques-neves-index}).  Being a central problem in the theory, we include this short section to prove an extension of  Main Theorem \ref{main.theorem} that, in light of Zhou's paper, is needed to obtain the precise characterization of the Morse index of min-max minimal hypersurfaces for generic metrics. This finishes the Morse-theoretic program (\cite{marques-icm},  \cite{marques-neves-cycles}, \cite{marques-neves-index},  \cite{neves-icm}) proposed by the  authors  for the area functional (see Theorem \ref{morse.program} below).

 Let $\mathcal{M}_C(g)$ denote the set of connected, closed, smooth, embedded minimal hypersurfaces (for the metric $g$) with both area and Morse index bounded from above by $C$.  The next  statement follows  from  Zhou's work:
\subsection{Multiplicity One Theorem}\label{zhou.theorem}(Zhou, \cite{zhou-multiplicity}) {\em Fix $k\in \N$ and consider  $(M^{n+1},g)$ a closed Riemannian manifold ,  $3\leq (n+1) \leq 7$. Suppose that   every element of $\mathcal{M}_{\omega_k(M,g)+k+1}(g)$ is nondegenerate. 

Then for every homotopy class $\Pi$ of $k$-sweepouts with ${\bf L}(\Pi)=\omega_k(M,g)$, there is a minimizing sequence $\{\Phi_j\}_j$ in $\Pi$ so that ${\bf C}(\{\Phi_j\}_j)$ contains a smooth, closed, embedded, multiplicity one, two-sided, minimal hypersurface $\Sigma$ with 
$$
\omega_k(M,g) = {\rm area}_g(\Sigma)\quad\mbox{and}\quad{\rm index}(\Sigma) \leq k.
$$
}
In Theorem A of \cite{zhou-multiplicity} it is assumed that $g$ is bumpy but an inspection of the proof shows that requiring every element of $\mathcal{M}_{\omega_k(M,g)+k+1}(g)$ to be nondegenerate suffices. The existence of the minimizing sequence $\{\Phi_j\}_j$ in $\Pi$ having $\Sigma\in {\bf C}(\{\Phi_j\}_j)$ is not specifically stated in \cite{zhou-multiplicity} but a look at Section 5 in \cite{zhou-multiplicity} shows that this is indeed the case. The delicate part of his theorem is the statement that $\Sigma$ has multiplicity one; the upper index bounds had been proven previously in \cite{marques-neves-index}.

An inspection of the  proof of  Theorem 1.3 yields the following improvement:

\subsection{Theorem}\label{deformation.theorem}
{\it Suppose $(M^{n+1},g)$ is a closed Riemannian manifold, $3\leq (n+1) \leq 7$, such that for some $k\in \mathbb{N}$, every element of $\mathcal{M}_{\omega_k(M,g)+k+1}(g)$ is nondegenerate.  Let $\Pi$ be a homotopy class of $k$-sweepouts with ${\bf L}(\Pi)=\omega_k(M,g)$. Suppose $\{\Phi_j\}$ is a minimizing sequence in $\Pi$ such that every embedded minimal cycle of ${\bf C}(\{\Phi_j\})$ with index less than or equal to $k$ has multiplicity one. Then
there exists an embedded minimal cycle  $\Sigma \in {\bf C}(\{\Phi_j\})$ (hence ${\rm area}(\Sigma)=\omega_k(M,g)$) with
$$
{\rm index}(\Sigma) = k.
$$
}

We will show how Theorem \ref{deformation.theorem} implies:

\subsection{Theorem}\label{improvement.theorem}
{\em Fix $k\in\N$ and assume Zhou's Multiplicity One Theorem (Theorem \ref{zhou.theorem}). 
 
For every bumpy metric $g$ on a closed manifold  $M^{n+1}$, $3\leq (n+1) \leq 7$, there exists a smooth, closed, embedded, multiplicity one, two-sided, minimal hypersurface $\Sigma_k$ such that
$$
\omega_k(M,g) = {\rm area}_g(\Sigma_k),\quad\mbox{and}\quad{\rm index}(\Sigma_k) = k.
$$
}

As a consequence,  when combined with  the Weyl Law for the Volume Spectrum proven in \cite{liokumovich-marques-neves} by Liokumovich and the authors, one finally obtains
the Morse-theoretic result predicted by the authors:

\subsection{Theorem}\label{morse.program}
{\it Let $g$ be a  $C^\infty$-generic (bumpy) metric  on a closed manifold  $M^{n+1}$, $3\leq (n+1) \leq 7$.  For each $k\in \mathbb{N}$, there exists a smooth, closed, embedded, multiplicity one, two-sided, minimal hypersurface $\Sigma_k$ such that
$$
\omega_k(M,g) = {\rm area}_g(\Sigma_k)\quad\mbox{and}\quad{\rm index}(\Sigma_k) = k
$$
and 
$$
\lim_{k\rightarrow \infty}\frac{{\rm area}_g(\Sigma_k)}{k^\frac{1}{n+1}}=a(n) {\rm vol}(M,g)^\frac{n}{n+1},
$$
where $a(n)>0$ is the dimensional constant in the Weyl law for the volume spectrum (\cite{liokumovich-marques-neves}).
}


\subsection{Proof of  Theorem \ref{improvement.theorem}}

 Sharp's  Compactness Theorem \cite{sharp} implies  the set $\mathcal{M}_C(g)$  is finite for every $C>0$,  if the metric $(M,g)$ is bumpy.

\subsection{Proposition}\label{rational.independence} {\em Let $g$ be a bumpy metric on a closed manifold $M^{n+1}$. Given $C>0$, there exists a sequence of Riemannian metrics $(g_i)_{i\in \mathbb{N}}$ converging to $g$ in the smooth topology so that for each $i\in \mathbb{N}$:
\begin{itemize}
\item every element of $\mathcal{M}_C(g_i)$ is $g_i$-nondegenerate;
\item and if 
$$
p_1 \cdot {\rm area}_{g_i}(\Sigma_1) + \cdots + p_N \cdot {\rm area}_{g_i}(\Sigma_N)=0,
$$
with $\{p_1, \dots, p_N\} \subset \mathbb{Z}$, $\{\Sigma_1, \dots, \Sigma_N\}\subset \mathcal{M}_C(g_i),$ and $\Sigma_k \neq \Sigma_l$ whenever $k\neq l$, then
$$
p_1 = \cdots = p_N = 0.
$$
\end{itemize}
 }

\begin{proof}
Since $g$ is bumpy, $\mathcal{M}_C(g)$ is finite. Let $\mathcal{M}_C(g)=\{S_1, \dots, S_q\}$, $q\in \mathbb{N}$, $S_k\neq S_l$ whenever $k\neq l$.

Recall that if $\tilde{g}=\exp(2\phi)g$, then the second fundamental form of $\Sigma$ with respect to $\tilde{g}$ is given by (Besse \cite{besse}, Section 1.163)
\begin{eqnarray*}\label{second.fundamental.form}
A_{\Sigma, \tilde{g}} =  A_{\Sigma,g} -  g \cdot (\nabla \phi)^\perp,
\end{eqnarray*}
where $(\nabla \phi)^\perp(x)$ is the component of $\nabla \phi$ normal to $T_x\Sigma$.
 
We can pick $p_l \in S_l \setminus (\cup_{k\neq l} S_k)$ for every $l=1, \dots, q$ (see the proof of Lemma 4 of 
\cite{marques-neves-song}). Let $\varepsilon>0$ be sufficiently small so that $B_\varepsilon(p_k)\cap B_\varepsilon(p_l) = \emptyset$ whenever $k\neq l$ and $B_\varepsilon(p_l) \cap (\cup_{k\neq l} S_k) = \emptyset$ for every $l=1, \dots, q$. We choose a nonnegative function $f_l\in C_c^\infty(B_\varepsilon(p_l))$, ${f_l}_{|S_l}\not \equiv 0$, such that $(\nabla_g f_l)(x) \in T_xS_l$ for every $x\in S_l$. Hence $S_l$ is still minimal with respect to the metric $\hat{g}(t_1,\dots,t_q) = \exp(2(t_1f_1+\cdots+t_qf_q))g$, for every $l=1, \dots, q$.

Let $(t_1^{(i)}, \dots, t_q^{(i)})\in (0,1]^q$ be a sequence converging to zero so that, by putting $g_i=\hat{g}(t_1^{(i)}, \dots, t_q^{(i)})$, we have that  the real numbers
$$
{\rm area}_{g_i}(S_1), \dots, {\rm area}_{g_i}(S_q)
$$
are linearly independent over $\mathbb{Q}$. Sharp's compactness theorem, together with the fact that $S_l$ is nondegenerate with respect to $g$ for every $l=1, \dots, q$, implies that for sufficiently large $i$ we have
$$
\mathcal{M}_C(g_i) =\{S_1, \dots, S_q\},
$$
and each $S_l$ is $g_i$-nondegenerate.
This finishes the proof of the lemma.
\end{proof}

A look at the proof of Proposition \ref{find.homotopy.class} leads to the improvement:

\subsection{Proposition}\label{find.homotopy.class.2}
{\it Suppose $(M^{n+1},g)$ is a closed Riemannian manifold, $3\leq (n+1) \leq 7$, such that for some $k\in \mathbb{N}$, every element of $\mathcal{M}_{\omega_k(M,g)+k+1}(g)$ is nondegenerate.  Then there exists a homotopy class $\Pi$ of $k$-sweepouts such that
$$
\omega_k(M,g)={\bf L}(\Pi).
$$
}


\begin{proof}[Proof of Theorem \ref{improvement.theorem}]
Fix $k\in \mathbb{N}$ and let $C=\omega_k(M,g)+k+2$. By Proposition \ref{rational.independence}, there exists a sequence of Riemannian metrics $(g_i)_{i\in \mathbb{N}}$ converging to $g$ in the smooth topology so that for each $i\in \mathbb{N}$:
\begin{itemize}
\item[(i)] every element of $\mathcal{M}_C(g_i)$ is $g_i$-nondegenerate;
\item[(ii)] and if 
$$
p_1 \cdot {\rm area}_{g_i}(\Sigma_1) + \cdots + p_N \cdot {\rm area}_{g_i}(\Sigma_N)=0,
$$
with $\{p_1, \dots, p_N\} \subset \mathbb{Z}$, $\{\Sigma_1, \dots, \Sigma_N\}\subset \mathcal{M}_C(g_i),$ and $\Sigma_k \neq \Sigma_l$ whenever $k\neq l$, then
$$
p_1 = \cdots = p_N = 0.
$$
\end{itemize}

For sufficiently large $i$, $\mathcal{M}_{\omega_k(M,g_i)+k+1}(g_i) \subset \mathcal{M}_C(g_i)$. Proposition \ref{find.homotopy.class.2} then implies the existence of a homotopy class $\Pi_i$ of $k$-sweepouts such that
$$
\omega_k(M,g_i)={\bf L}(\Pi_i).
$$

From Theorem \ref{zhou.theorem} (Zhou's Multiplicity One Theorem), we obtain  a minimizing sequence $\{\Phi_j^{(i)}\}_j$  in $\Pi_i$ such that the critical set ${\bf C}(\{\Phi_j^{(i)}\}_j)$ contains a multiplicity one, two-sided, embedded $g_i$-minimal cycle 
$$
V_i= \Sigma_1^{(i)} + \cdots +\Sigma_{N_i}^{(i)},
$$
$\{\Sigma_1^{(i)}, \dots, \Sigma_N^{(i)}\}\subset \mathcal{M}_C(g_i),$ and $\Sigma_k^{(i)} \neq \Sigma_l^{(i)}$ whenever $k\neq l$,
 and
 $$
 {\rm index}(V_i) = {\rm index}(\Sigma_1^{(i)}) + \cdots +  {\rm index}(\Sigma_{N_i}^{(i)}) \leq k.
 $$
 
 Now let 
 $$
 \tilde{V}_i = m_1^{(i)} \cdot \tilde\Sigma_1^{(i)} + \cdots +m_{Q_i}^{(i)}\cdot \tilde\Sigma_{Q_i}^{(i)}, \, \, \{m_1^{(i)}, \dots, m_{Q_i}^{(i)}\}\subset \mathbb{N}
 $$
 be any embedded minimal cycle in ${\bf C}(\{\Phi_j^{(i)}\}_j)$ with index less than or equal to $k$. Then $\{\tilde\Sigma_1^{(i)}, \dots, \tilde\Sigma_{Q_i}^{(i)}\} \subset \mathcal{M}_C(g_i)$. Since any element in ${\bf C}(\{\Phi_j^{(i)}\}_j)$ has area equal to $\omega_k(M,g_i)$, we have 
 $$
 {\rm area}(V_i)={\rm area}(\tilde{V}_i).
 $$
 By property (ii) above, we must have $\tilde{V}_i=V_i$.
 
 This implies that the assumptions of Theorem \ref{deformation.theorem} are satisfied for the metric $g_i$ and the minimizing sequence $\{\Phi_j^{(i)}\}_j$. Hence there exists a multiplicity one embedded minimal cycle  $\Sigma_k^{(i)} \in {\bf C}(\{\Phi_j^{(i)}\}_j)$ with
$$
{\rm index}(\Sigma_k^{(i)}) = k.
$$
Necessarily $\Sigma_k^{(i)}=V_i$, hence it is two-sided. Note that ${\rm area}_{g_i}(\Sigma_k^{(i)})=\omega_k(M,g_i)$. By Sharp's compactness theorem, $\{\Sigma_k^{(i)}\}$ converges subsequentially to a $g$-minimal hypersurface $\Sigma_k$ which must have multiplicity one and index equal to $k$ because $g$ is bumpy.  This finishes the proof of the theorem.
\end{proof}

\bibliographystyle{amsbook}

\end{document}